\documentclass[11pt]{amsart}
\pdfoutput=1
\usepackage[utf8]{inputenc}
\usepackage{amsmath, amssymb}
\usepackage[english]{babel}
\usepackage{graphicx}
\usepackage{amsthm}
\usepackage{array}
\usepackage{amsthm}
\usepackage{mathtools}
\usepackage{thmtools}
\usepackage{lpic}
\usepackage{comment}
\usepackage[margin=1 in]{geometry}
\allowdisplaybreaks
\usepackage{enumitem}
\usepackage{subcaption}
\usepackage[dvipsnames]{xcolor}
\usepackage{hyperref}
\definecolor{candyapplered}{rgb}{1.0, 0.03, 0.0}
\definecolor{mediumblue}{rgb}{0.0, 0.0, 0.8}
\hypersetup{
pdfauthor={Tarik Aougab, Jonah Gaster},
pdftitle={},
bookmarksnumbered,
colorlinks=true,
linkcolor=mediumblue,
citecolor=candyapplered,
urlcolor=blue}

\declaretheorem[numberwithin=section]{theorem}
\declaretheorem[numbered=no,name=Main Theorem]{main}
\declaretheorem[sibling=theorem, style=definition]{definition}
\declaretheorem[sibling=theorem]{lemma}
\declaretheorem[sibling=theorem, style=remark]{remark}
\declaretheorem[sibling=theorem]{proposition}

\declaretheorem[style=remark]{question}
\declaretheorem[numbered=no,style=remark]{Question}

\numberwithin{equation}{section}

\newcommand{\bZ}{\mathbb{Z}}

\newcommand{\bH}{\mathbb{H}}

\newcommand{\cC}{\mathcal{C}}
\newcommand{\cP}{\mathcal{P}}
\newcommand{\cA}{\mathcal{A}}
 
\newcommand{\cN}{\mathcal{N}}

\newcommand{\cF}{\mathcal{F}}

\newcommand{\cM}{\mathcal{M}}

\newcommand{\Area}{\mathrm{Area}}

\newcommand{\Stab}{\mathrm{Stab}}

\title{From arcs to curves: Quadratic growth of 1-systems}
\author{Tarik Aougab}
\address{Department of Mathematics \\ Haverford College \\ Haverford, PA 19041}
\email{taougab@haverford.edu}

\author{Jonah Gaster}
\address{Department of mathematical sciences \\ University of Wisconsin-Milwaukee \\ Milwaukee, WI 53211} 
\email{gaster@uwm.edu}

\date{November 30, 2026}

\begin{document}

\begin{abstract}
We show that a collection of simple closed curves pairwise intersecting at most once on an orientable surface of Euler characteristic $\chi$ has at most $2016|\chi|^2+338|\chi|$ curves.
This resolves a thirty-year old problem, up to multiplicative constants. Inspired by the work of Przytycki in the setting of arcs, we introduce the concepts of \emph{tulips}, \emph{flowers}, and \emph{stem systems} in order to account for how certain polygons built from pairs of curves in the collection distribute area over the surface. 
\end{abstract}

\maketitle

\section{introduction}

Let $\Gamma$ be a collection of homotopy classes of simple closed curves on an orientable surface $S$ of Euler characeristic $\chi$. We call $\Gamma$ a $k$-system provided the geometric intersection number  $\iota(\gamma,\gamma')$ is bounded by $k$ for every $\gamma,\gamma'\in\Gamma$. Rough estimates for the size of $k$-systems on surfaces were obtained first by Juvan-Malni\c{c}-Mohar \cite{JMM}. 
Determining the growth rate for the case $k=1$, as a function of $\chi$, is surprisingly difficult. 
Attention on this question was 
spurred by Farb-Leininger, 
and more recently by the \emph{K3 Problem List} \cite[Problem 2.12(b)]{k3}:

\begin{Question}
How large may a $1$-system of curves be on $S$?
\end{Question}

Constructions of size quadratic in $|\chi|$ were known to Farb-Leininger, but upper bounds have been much harder to come by. Malestein-Rivin-Theran produced exponential upper bounds \cite{MRT}; Przytycki obtained cubic upper bounds \cite{Przytycki}; and Aougab-Biringer-Gaster obtained sub-cubic upper bounds \cite{ABG}. Finally, the best upper bounds prior to this work are due to Greene, who used a probabilistic argument to produce a bound of the form $O(|\chi|^2\log |\chi|)$ \cite{Greene1}.

In this paper, we show:

\begin{main}
\label{main thm}
If $\Gamma$ is a 1-system on the orientable surface $S$ with Euler characteristic $\chi<0$, 
then the size $|\Gamma|$ is at most 
$84 |\chi|(24|\chi|+4)+\frac32|\chi|$.
\end{main}

As a corollary, the size of the largest $1$-system on $S$ grows quadratically in $\chi$. 
The last term above can be more precisely replaced by the size of a pants decomposition on $S$. 
Of course, this is irrelevant because the first term dominates, and in any case the coefficient of $|\chi|^2$ is surely not optimal (see \S\ref{sec:questions}). 
We strongly believe 
the strategies we use to prove 
our main theorem will produce upper bounds for $k$-systems on the order of 
$|\chi|^{k+1}$. 
This will be explored in a 
follow-up paper.

\subsection{Strategy}
Przytycki's approach to a cubic bound for the size of a 1-system relied on 
the analogous count for \emph{arcs} on a punctured surface $S$. 
Incredibly, Przytycki computed the maximum size of a 1-system of arcs on a punctured surface exactly. Since Przytycki's remarkable paper, all subsequent improvements to the 1-systems bound relied on his work as a black box. Our central idea is to open up this box and modify its components to achieve a similar level of control for curves on surfaces. We will briefly review Przytycki's argument in \S\ref{sec:definitions}. 

For the reader familiar with Przytycki's work, our strategy can be summed up as follows: choose a maximal multicurve $\cC\subset \Gamma$, and consider what happens when the curves in $\cC$ are pinched to cusps. 
The remaining curves turn into arcs that give rise to \emph{nibs}, as in Przytycki. However, a fundamental stumbling block is that a curve might pass through many curves in $\cC$, and so may turn into many arcs, making estimation of $|\Gamma\setminus\cC|$ in terms of nibs impossible. 
The solution is to not pinch $\cC$ \emph{completely}, but to stop just short. 
Elements of $\Gamma\setminus \cC$ may pass through many curves in $\cC$, but it turns out that by leveraging the hyperbolic geometry of the picture-- in particular, by paying careful attention to when pairs of curves fellow-travel and when they diverge-- there is a way to pick out `base curves' in $\cC$ for each of the curves in $\Gamma\setminus \cC$, from which one may observe what seem to be \textit{approximate} nibs and \textit{approximate} slits. In our present terminology these become \emph{tulips}, \emph{flowers}, and \emph{stems}. This terminology will be explained in \S\ref{sec:definitions} and \S\ref{sec:maximal}.

\begin{theorem}
\label{thm:stems bound}
With $\Gamma$ and $\cC$ as above, we have 
\begin{equation}
\label{eq:main inequality}
|\Gamma|\;\le 
\;|\cC|\;+\;
84|\chi|\,\sigma(S)~,
\end{equation}
where $\sigma(S)$ is the size of the largest stem system on $S$, a topological constant depending only on $S$.
\end{theorem}

We refer to \S\ref{sec:stems} for the precise meaning of a stem system. 
Roughly speaking, it is a collection of arcs (i.e.~\emph{stems}) from a point $p\in S$ to a nonempty multicurve $\cC\subset S$ 
with the property that 
no stem is a prefix of another 
(up to homotopies of $(I,\partial I)\to (S,\{p\}\cup\cC)$), 
and one more technical condition, which the reader can think of for now as almost saying that distinct stems are disjoint--the latter needs to be interpreted carefully 
with respect to $\cC$.

In context, the multicurve $\cC$ is drawn from the 1-system $\Gamma$,
and 
a stem system is what one sees when a point is in many different 
approximate nibs (that is, \emph{tulips}) satisfying a certain maximality condition. 
However, note that in the wake of Theorem~\ref{thm:stems bound} the $1$-system $\Gamma$ is irrelevant: $\sigma(S)$ is obtained by maximizing over \emph{any} stem system attached to \emph{any} multicurve $\cC$. 

\begin{proposition}
\label{prop:sigmag}
We have $2|\chi|+2\le\sigma(S)\le 24|\chi|+4$.
\end{proposition}

The main theorem follows immediately from Theorem~\ref{thm:stems bound} and Proposition~\ref{prop:sigmag}. We will describe the proof of Theorem~\ref{thm:stems bound} in the following section while motivating some of the terminology, and the proof of Proposition~\ref{prop:sigmag} can be found in \S\ref{sec:stems}.

\begin{remark}
The lower bound for $\sigma(S)$ above is not relevant to our main theorem, but it demonstrates linear growth of $\sigma(S)$ in $|\chi|$, and indicates limitations to improving the upper bounds we obtain for $|\Gamma|$. 
We are unable to determine where $\sigma(S)$ lies between the two linear bounds, though we are interested in having a precise calculation. 
While it is not hard to determine $\sigma(S_1)=2$, where $S_g$ is the closed orientable surface of genus $g$, we are even unsure of $\sigma(S_2)$.

We suspect $
\sigma(S)$ actually lies closer to the given lower bound. In fact, an earlier version of this manuscript included an upper bound of $6|\chi|+6$, with a claimed proof that attempted to discard hyperbolic geometry entriely. 
However, the proof had holes spotted by a careful referee. 
Here we present a corrected proof that is more transparent, and philosophically more in line with the rest of the paper, but at the expense of a larger slope in the upper bound.
\end{remark}

\subsection{Outline.} In \S\ref{sec:definitions}, we summarize Przytycki's approach in the setting of arcs and introduce the concepts of \emph{tulips}, \emph{flowers}, \emph{stems}, and \emph{stalks} that we will need throughout the paper. 

One complication that arises in our setting is that, in some sense, there are too many tulips. For our strategy to work, we have to focus on a subset of them which are ``maximal'' -- we outline this in \S\ref{sec:maximal}. This however introduces another complication, because we do not know that every curve in $\Gamma$ arises as a side of a maximal tulip, and so it becomes difficult to deduce a bound on $|\Gamma|$ by only studying the maximal tulips. The solution is to categorize tulips into \emph{types} which, loosely speaking, keep track of how the sides of a tulip twist around a certain curves of $\cC$. At the end of \S\ref{sec:maximal}, we show that every curve is a main curve of a tulip that is maximal \emph{amongst its type} (i.e.~\emph{type-maximal}); combining this with a bound on the number of possible types allows us to argue type-by-type.  

In \S\ref{sec:technical lemma}, we 
demonstrate that a collection of type-maximal tulips gives rise to a stem system. 
In \S\ref{sec:stems}, we bound the size of a stem system and show that it can be at most linear in $|\chi|$.  
Finally, in \S\ref{sec:theorem proofs}, we put everything together to prove Theroem~\ref{thm:stems bound}, and in \S\ref{sec:questions} we gather some remaining questions about the structure and size of $1$-systems.

\subsection{Acknowledgements.} 
The authors would like to thank everyone who has listened to either one of them talk about this problem over the last 10 years. We especially extend our gratitude to Ian Biringer, Josh Greene, and Piotr Przytycki, with whom many helpful and inspiring conversations were enjoyed, and to a referee, whose critically useful suggestions have improved this work.

\section{Tulips, flowers, stalks, and stems}
\label{sec:definitions}

\subsection{Przytycki's approach}
\label{sec:narration}
Suppose that $\cA$ is a set of simple proper geodesic arcs on a cusped hyperbolic surface $X$. Around a cusp, ends of elements of $\cA$ are cyclically ordered, and consecutive elements give rise to an immersion $\Delta\looparrowright X$ of an ideal triangle $\Delta\subset \bH$; this is called a \emph{nib} for $\cA\subset X$. Moreover, a nib has a preferred corner of the ideal triangle $\Delta$, and so each point $q\in \Delta$ gives rise to a proper geodesic ray from $q$ to $\partial_\infty\bH$, whose image under the nib map is called a \emph{slit}. 

Przytycki's bounds are obtained by considering fibers of the local isometry
\[
 \nu:\bigsqcup_{nibs} \Delta \looparrowright X~.
\]
For $p\in X$, the fiber $\nu^{-1}(p)$ 
gives rise to a collection of slits from $p$ to the cusps of $X$. When elements of $\cA$ pairwise intersect at most once, Przytycki's count for $|\cA|$ is achieved with two observations:
\begin{enumerate} 
\item
\label{it:P1}
The size $|\nu^{-1}(p)|$ is at most $2(|\chi|+1)$, and the number of nibs is at most $4|\chi|(|\chi|+1)$.
\item 
\label{it:P2}
The number of arcs $|\cA|$ is equal to one-half the number of nibs.
\end{enumerate}

In Przytycki's work, the deep realization is certainly \eqref{it:P1}, which is a consequence of the two beautiful facts that slits are embedded, and, provided $\cA$ is a 1-system, are pairwise disjoint-- the nibs bound comes from considering hyperbolic area. 
Item \eqref{it:P2} is an obvious consequence of the construction, and provides the claimed bound for $|\cA|$.
Przytycki's cubic bound for 1-systems follows quickly as well: let $\cC\subset \Gamma$ be a maximal multicurve. The subset $\Gamma_c=\{\gamma\in \Gamma:\iota(\gamma,c)\ne 0\}$ maps at worst 2-to-1 to a 1-system of arcs on the surface obtained by pinching $c$ to a pair of cusps, so $|\Gamma_c|$ is bounded quadratically in $|\chi|$. Because $|\cC|\le \frac32|\chi|$, the size of $\Gamma=\bigcup_{c\in\cC}\; \Gamma_c$ is $\lesssim |\chi|^3$.

In order to adapt Przytycki's work for closed curves, we must first construct the `nibs', and then analyze fibers of the map from `nibs' to the surface. These two steps pull against each other: in order to account for all the curves, we would like to flexibly allow many `nibs'. However, if there are too many `nibs', the fibers will grow too large, complicating part \eqref{it:P1}. Conversely, if `nibs' are defined restrictively so that the fibers are small, it will be difficult to ensure we are capturing enough of the curves to obtain a good estimate for $|\Gamma|$, interfering with part~\eqref{it:P2}. 

Before explaining the terminology, we highlight modification of Przytycki's strategy:

\begin{theorem}
\label{thm:P1}
The number of tulips that are type-maximal is at most $84|\chi|\sigma(S)$, where $\sigma(S)$ is the size of the largest stem system on $S$.
\end{theorem}

\begin{theorem}
\label{thm:P2}
The number of curves in $\Gamma\setminus\cC$ is at most the number of type-maximal tulips.   
\end{theorem}

Theorems~\ref{thm:P1} and \ref{thm:P2} together prove Theorem~\ref{thm:stems bound}.

We turn now to introduce the construction of tulips and related terminology.

\subsection{Choosing a hyperbolic metric}
Throughout the following, $S$ is a closed orientable surface of genus $g\ge 2$, and $\Gamma$ is a 1-system on $S$. 

Let $\cC$ be any maximal multicurve in $\Gamma$, and let $\epsilon>0$.
As $\cC$ is maximal, each curve in $\Gamma\setminus \cC$ intersects a curve in $\cC$.
It will be useful as well to have a pants decomposition on $S$, so we arbitrarily enlarge $\cC$ to a pants decomposition $\cP_\cC$; of course, the simpler scenario $\cC=\cP_\cC$ is worthwhile to keep in mind.
We start by choosing a hyperbolic metric on $S$ so that the lengths of the curves in $\cP_\cC$ are smaller than $\epsilon$, 
an immediate consequence of the nature of Fenchel-Nielsen coordinates \cite[Ch.~3.6]{Buser}.

\begin{lemma} \label{lem:metric}
For every $\epsilon>0$ there is a hyperbolic structure on $S$ so that each curve in $\cP_\cC$ has length at most $\epsilon$. 
\end{lemma}

Using this hyperbolic structure, we identify the universal covering with $\pi:\bH\to S$, and we replace all curves in $\Gamma$ with their geodesic representatives. Observe that we may endow the geodesics in $\Gamma$ with orientations arbitrarily, according to convenience.

\subsection{Cyclic orders on annuli}
Let $c\in\cC$.
Choose an orientation for $c$, and a lift $\tilde{c}\subset\bH$, and let $\Gamma_{\widetilde{c}}$ indicate the collection of lifts of curves in $\Gamma_c$ that intersect $\tilde{c}$. Endow all geodesics in $\Gamma_{\widetilde{c}}$ with orientations so they interesect $\tilde{c}$ positively, and observe that the endpoints of the curves in $\Gamma_{\widetilde{c}}$ are now linearly ordered by $\partial_\infty\bH$. Because this linear order is preserved by the cyclic subgroup of $\pi_1 S$ fixing $\tilde{c}$, we obtain a cyclic order for the curves in $\Gamma_c$, evidently independent of the choice of $\tilde{c}$.
Because $\Gamma_c$ is finite, we may immediately conclude: 

\begin{lemma}
\label{lem:discrete}
The linear order $\le$ on $\Gamma_{\widetilde{c}}$ is discrete.
\end{lemma}

\begin{remark}
The reader should notice that the cyclic order we obtain on $\Gamma_c$ is almost certainly different than the cyclic order on the intersection points $\{\gamma\cap c:\gamma\in\Gamma_c\}$ induced by $c$. In particular, it is possible that curves from $\Gamma$ form a triangle with $c$, in such a way that the cyclic order on $\Gamma_c$ induced by the order in the universal cover of the endpoints of their lifts along $\partial_\infty\bH$ is different than their order along $c$.
Philosophically speaking, it is especially pleasing that the order we use above is evidently independent of the chosen hyperbolic metric; indeed the order is determined by the positions of points on $\partial_\infty\Gamma\approx S^1$. 
\end{remark}

\subsection{The floral invariants}
Consider the dual tree $G(\cP_\cC)$ to $\pi^{-1}\cP_\cC$, with a vertex for each complementary component of $\pi^{-1}\cP_\cC\subset\bH$, and an edge for a lift of an element of $\cP_\cC$ between neighboring components of $\bH\setminus\pi^{-1}\cP_\cC$. The choice of lift $\tilde{c}$ for $c\in\cC$ determines an edge of $G(\cP_\cC)$, and each of the elements $\alpha\in\Gamma_{\widetilde{c}}$ determine a trajectory in $G(\cP_\cC)$ that goes through this edge.
Let us indicate the sequence of edges of $G(\cP_\cC)$ visited by $\alpha\in \Gamma_{\widetilde{c}}$ by the bi-infinite sequence $E(\alpha)=(a_j)$.

\begin{figure}[ht]
\begin{lpic}{pics/nibConstructionAgain(,7cm)}
\Large
\lbl[]{40,56;$\tilde{c}$}
\lbl[]{60,150;$\widetilde{c_1}$}
\lbl[]{108,155;$\widetilde{c_2}$}
\lbl[]{130,32;$\widetilde{c_-}$}
\lbl[]{35,83;$\widetilde{c_0}$}
\end{lpic}
\caption{The tulip between $\alpha$ and $\beta$ using the dual graph $G(\cP_\cC)$. Lifts of curves in $\cP_\cC$ are distinguished as to whether they are in $\Gamma$ (full stroke) or not (dotted).}
\label{fig:tulip}
\end{figure}

Suppose that $\widetilde{\beta}$ immediately follows $\widetilde{\alpha}$ in the linear order on $\Gamma_{\widetilde{c}}$, with trajectories $E(\widetilde{\alpha})=(a_j)$ and $E(\widetilde{\beta})=(b_j)$. 
Let $\lambda_0$ indicate the maximal shared subsequence of $(a_j)$ and $(b_j)$, nonempty because $\widetilde{\alpha},\widetilde{\beta}\in\Gamma_{\widetilde{c}}$ and finite because endpoints of lifts of distinct closed geodesics are disjoint \cite[p.~32]{Farb-Margalit}. Finally, we may consider the suffix $\lambda\subset \lambda_0$ with the property that the first entry of $\lambda$ is the first coordinate of $\lambda_0$ that is an element of $\cC$; observe that such an entry exists again because $\widetilde{\alpha},\widetilde{\beta}\in\Gamma_{\widetilde{c}}$.
After shifting $E(\widetilde{\alpha})=(a_j)$ and $E(\widetilde{\beta})=(b_j)$ simultaneously, we may assume that the entry $a_0=b_0$ is the \emph{last} entry of $\lambda$.  Observe that $\lambda$ now takes the form 
\[
\lambda=(a_{-r},\ldots,a_0)=(b_{-r},\ldots,b_0)~,
\]
where $r$ is the length of $\lambda$, and note that $a_1\ne b_1$.

We fix some notation for certain lifts in $\pi^{-1}\cP_\cC$:
The edge of $G(\cP_\cC)$ at $a_0=b_0$ is given by $\widetilde{c_0}$, the edges of $G(\cP_\cC)$ at $a_1$ and $b_1$ are given by $\widetilde{c_1}$ and $\widetilde{c_2}$, respectively, and the edge of $G(\cP_\cC)$ at $a_{-r}=b_{-r}$ is given by $\widetilde{c_-}$.
We stress: 
\begin{equation}
\tag{$\ast$}
\label{eq:assumption}
\text{\emph{
By construction, $\widetilde{c_-}$ is a lift of an element in $\cC$~.
}}
\end{equation}
See Figure~\ref{fig:tulip} for a picture of this setting.

\begin{definition}[Tulips, flowers, and stalks]
\label{def:nibs etc}
Let $\alpha$ and $\beta$ be consecutive elements in the cyclic order on $\Gamma_c$, and suppose lifts  $\widetilde{\alpha}$, $\widetilde{\beta}$, $\widetilde{c}$, $\widetilde{c_0}$, $\widetilde{c_1}$, $\widetilde{c_2}$, $\widetilde{c_-}$ are chosen as above.

The \emph{tulip} between $\alpha$ and $\beta$ is the local isometry $\pi|_N:N\to S$, 
where $N$ is the quadrilateral given by the convex hull of the four points $\widetilde{\alpha}\cap\widetilde{c_-}$, $\widetilde{\beta}\cap\widetilde{c_-}$,
$\widetilde{\alpha}\cap\widetilde{c_1}$, and
$\widetilde{\beta}\cap\widetilde{c_2}$.

The \emph{flower} of the tulip between $\alpha$ and $\beta$ is the 
local isometry $\pi|_F:F\to S$, where $F$ is the quadrilateral given by 
the intersection of $N$ with the halfspace bounded by $\widetilde{c_0}$ and containing $\widetilde{c_i}$.

The \emph{stalk} of the tulip between $\alpha$ and $\beta$ is the 
local isometry $\pi|_Q:Q\to S$, where $Q$ is the quadrilateral given by 
the intersection of $N$ with the halfspace bounded by $\widetilde{c_0}$ and containing $\widetilde{c_-}$. In the event that $\widetilde{c_0} = \widetilde{c_-}$, we declare the stalk to be empty. 

See Figure~\ref{fig:stem} for a schematic of this construction.
\end{definition}

As $\epsilon$ goes to $0$, either side of a collar neighborhood of a curve in $\cP_\cC$ geometrically limits to a neighborhood of a cusp. Simultaneously, the quadrilaterals given by the tulips, or their flowers, limit geometrically to ideal triangles. 
On the other hand,  were we to actually complete this geometric limit, the stalk of the tulip collapses to potentially several infinite simple geodesic arcs between cusps. The literal translation to Przytycki's setting is somewhat lost, as the slits get cut off. Thus our strategy is not to complete the geometric limit.
Nonetheless, there is enough here to carry through much of Przytycki's argument, starting with \cite[Lemma 2.5]{Przytycki}.

We will call the curves $\alpha$ and $\beta$ the \emph{main curves} of the tulip $N$; we refer to $c_-\in\cC$ as the \emph{base} of $N$, and $c_0\in \cP_\cC$ as the  \emph{base of the flower}.
The two vertices $\widetilde{\alpha}\cap \widetilde{c_1}$ and $\widetilde{\beta}\cap \widetilde{c_2}$ of $N$ are called the \emph{left} and \emph{right tops} of the tulip, respectively.
Two of the sides of $N$ are homotopic \emph{rel $\widetilde\cC$} to (lifts of) subarcs of the main curves; these are the \emph{main sides} of the tulip. 
Observe that $\widetilde{\alpha}$ and $\widetilde{\beta}$ may cross $N$ diagonally, i.e.~the main sides may not be subarcs of the main curves. 

\begin{remark}
\label{rem:notation}
Just as in Przytycki, we allow $\alpha=\beta$ in Definition~\ref{def:nibs etc} above, and we note that the construction is well-defined and independent of lifts chosen.
When clear from context, we may also refer to the lifts $\widetilde{\alpha}$ and $\widetilde{\beta}$ as the main curves of a tulip, the lifts $\widetilde{c_-}$ and $\widetilde{c_0}$ as the base and base of the flower, respectively, and the lifts $\widetilde{c_1}$ and $\widetilde{c_2}$ as the tops.
\end{remark}

\begin{remark}
\label{rem:pants decomp}
The difference in the roles of $\cC$ and $\cP_\cC$ may seem slightly mysterious for the moment. 
The reason why $\cP_\cC$ is necessary in our argument will only become clear in \S\ref{sec:maximal}. The reader should note, for instance, that because $\cC\subset \Gamma$, a main curve of a tulip intersects the base exactly once; on the other hand, if $\cP_\cC\ne\cC$ it may intersect the base of the flower arbitrarily many times.
Nonetheless, the main thrust of the argument can be appreciated in the case that $\cC=\cP_\cC$, which the reader is encouraged to assume if simplification is needed.
\end{remark}

If the main sides of a tulip project to sub-arcs of the main curves, then it is obvious that they (that is, their images) do not self-cross. 
This is true in any case:

\begin{lemma}
\label{lem:main sides}
The main sides of a tulip embed under $\pi|_N:N\to S$.
\end{lemma}

\begin{proof}
If a main side $\gamma$ of a tulip did not embed under $\pi$, there would be an element $g\in\pi_1S$ so that $\gamma$ and $g\cdot \gamma$ intersect transversely. The endpoints of $\gamma$ lie on certain lifts $\widetilde{c_i}$ and $\widetilde{c_-}$ from $\pi^{-1}\cC$, and the intersection $\gamma\cap g\cdot\gamma\ne \emptyset$ implies that the pairs $(\widetilde{c_i},\widetilde{c_-})$ and $(g\cdot\widetilde{c_i},g\cdot\widetilde{c_-})$ separate each other at infinity. This implies that $\widetilde{\alpha}$ and $g\cdot\widetilde{\alpha}$ intersect transversely, where $\alpha \in \Gamma$ is the corresponding main curve, a contradiction because $\alpha$ is simple. 
\end{proof}

\subsection{Stems}
In Przytycki's bounds, fibers of the map from nibs to the surface are bounded by the size of a collection of slits at a point $p\in S$. In our setting, nibs are replaced by tulips, with potentially very long stalks, and slits are replaced with \emph{stems at $\widetilde{p}$}. 

\begin{definition}[Stems at $\widetilde{p}$]
Suppose that $\widetilde{p}$ is in (the domain of) a tulip. The \emph{stem at $\widetilde{p}$} is the geodesic from $\widetilde{p}$ to the midpoint of the base of the tulip. See Figure~\ref{fig:stem}.
\end{definition}

\begin{figure}
\centering
\begin{lpic}{pics/nibTallAgain(,9cm)}
\Large
\lbl[]{100,135;$F$}
\lbl[]{100,58;$Q$}
\end{lpic}
\caption{A tulip, its flower $F$ and stalk $Q$, and a stem from a point in the flower.}
\label{fig:stem}
\end{figure}

One of the beautiful realizations in Przytycki's work is the fact that every slit in a nib embeds on the surface \cite[Lemma 2.5]{Przytycki}. 
One can almost copy-and-paste Przytycki's elegant proof, replacing the \emph{tips} of nibs with the \emph{bases} of tulips, which amounts to enlarging a point on $\partial_\infty\bH$ to an interval. 
However, the setting of tulips leaves a bit of flexibility not present in Przytycki's setting: distinct intervals can nest, whereas distinct points do not.
We require a small assumption to complete Przytycki's proof, which in any case is satisfied by `most' stems.

\subsection{The inner flower} 
Each flower of each tulip limits to an ideal hyperbolic triangle of area $\pi$ as $\epsilon\to0$. In fact, there is an especially well-behaved subset of the flower which does the same.

\begin{figure}[ht]
\includegraphics[width=7cm]{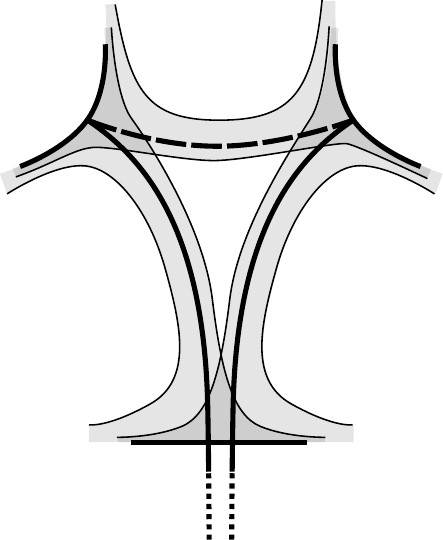}
\caption{The inner flower of a tulip is the complement in the flower of the shaded convex hulls pictured, formed by neighboring pairs of lifts of curves in $\cP_\cC$.}
\label{fig:innerFlower}
\end{figure}

\begin{definition}
[Inner flower]
\label{def:inner flower}
The \emph{inner flower} of a tulip is the geodesic triangle in the flower $F$ 
given by the complement of the three ideal quadrilaterals formed by the convex hulls of $\widetilde{c_0}\cup\widetilde{c_1}$, $\widetilde{c_0}\cup\widetilde{c_2}$, and $\widetilde{c_1}\cup\widetilde{c_2}$, respectively.
See Figure~\ref{fig:innerFlower}.
\end{definition}

In fact, the inner flower can be empty. In context, the elements of $\cP_\cC$ have length at most $\epsilon$, so by the Collar Lemma \cite{Farb-Margalit} the pairwise distances between the lifts $\widetilde{c_0}$, $\widetilde{c_1}$, and $\widetilde{c_2}$ approach $\infty$ as $\epsilon$ goes to $0$. Hence the flower of the tulip limits geometrically to an ideal triangle. Crucially for our bounds, the inner flower does as well, and we conclude:

\begin{lemma}
\label{lem:inner flower area}
As $\epsilon\to0$, the area of the inner flower approaches $\pi$.\qed
\end{lemma} 

In all that follows, we assume that $\epsilon$ is small enough so that all inner flowers have positive area.

The inner flower is useful for several reasons, one of which is that it provides a proper setting to generalize embeddedness of slits. 

\begin{remark}
\label{rem:tulip construction assumption}
Let us point out one essential ingredient below, without which the proof breaks down: By construction \eqref{eq:assumption}, the main curves of a tulip intersect its base only once on the surface. 
It could be interesting to pursue the construction of tulips by beginning with an \emph{arbitrary} pants decomposition, and dispensing with \eqref{eq:assumption}.
Without this assumption, however, we are unsure whether Proposition~\ref{prop:stems embed} remains true.
\end{remark}

\begin{proposition}
\label{prop:stems embed}
Stems from points in the inner flower embed under the covering map $\pi:\bH\to S$.
\end{proposition}

Recall that $\partial_\infty \bH\approx S^1$, which we orient naturally; given points $x,x'\in \partial_\infty \bH$ we indicate the oriented sub-interval of points from $x$ to $x'$ by $[x,x']$.
Given an oriented interval $I$ in $\partial_\infty\bH$, we indicate its endpoints by $I^\pm$, respectively.
We adapt \cite[Lemma 2.5]{Przytycki}:

\begin{proof}
Let $N\to S$ be a tulip, regarded as a subset of $\bH$, with base $\widetilde{c_-}$, tops $\widetilde{c_1}$ and $\widetilde{c_2}$, and a point $\widetilde{p}\in N$.
The stem at $\widetilde{p}$ is a complete geodesic segment from $\widetilde{p}$ to $\widetilde{c_-}$; it is contained in a complete geodesic segment $\gamma$ from the top arc of the inner flower to $\widetilde{c_-}$. 
We will show that $\gamma\to S$ is an embedding.
Towards contradiction, suppose there is $g\in\pi_1S$ so that $\gamma$ intersects $g(\gamma)$ transversely.

Each complete geodesic $\alpha$ in $\bH\setminus \{\widetilde{p}\}$ induces an oriented sub-interval of $\partial_\infty \bH$ by taking the closure of the side of the complement of $\alpha$ not containing $\widetilde{p}$, which we indicate by $[\alpha]$. Without loss of generality we may assume that $[\widetilde{c_1}]$, $[\widetilde{c_-}]$, $[\widetilde{c_2}]$ are positively oriented; see the left side of Figure~\ref{fig:Przytycki}.

First we consider $g(\widetilde{c_-})$, a lift in $\pi^{-1}\cP_\cC$.
If $g(\widetilde{c_-})=\widetilde{c_-}$, then $g(\gamma)$ and $\gamma$ would be disjoint. 
Hence $g(\widetilde{c_-})$ and $\widetilde{c_-}$ are disjoint. 
We claim that the intervals $[g(\widetilde{c_-})]$ and $[\widetilde{c_-}]$ are disjoint in $\partial_\infty \bH$.
Indeed, the alternative implies that either $[g(\widetilde{c_-})]\subset[\widetilde{c_-}]$ or $[\widetilde{c_-}]\subset[g(\widetilde{c_-})]$. In either case, we will argue that 
each of the projections of the main curves to the surface would intersect the projection of $\widetilde{c_-}$ at least twice, contradicting our construction \eqref{eq:assumption}. 

Indeed, suppose that $[\widetilde{c_-}]\subset[g(\widetilde{c_-})]$, so that $\gamma$ crosses $g(\widetilde{c_-})$ on its way to $\widetilde{p}$. 
Observe that no tops of flowers are separated by lifts in $\pi^{-1}\cP_\cC$; hence both main curves of $N$ intersect both $\widetilde{c_-}$ and $g(\widetilde{c_-})$.
By Lemma~\ref{lem:main sides}, these two projections are not $\pi_1S$-equivalent, so each of the main curves of $N$ would intersect the projection of $\widetilde{c_-}$ twice.
For the alternative $[g(\widetilde{c_-})]\subset [\widetilde{c_-}]$, the same argument applies to the main curves of $g(N)$.

\begin{figure}[ht]
\begin{minipage}{.45\linewidth}
\centering
\vspace{.8cm}
\begin{lpic}{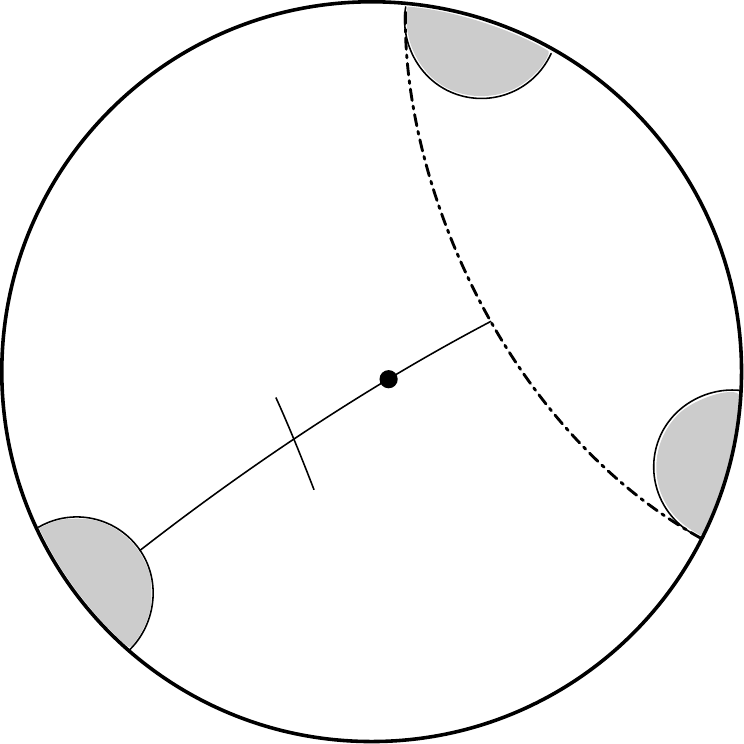(7cm)}
\Large
    \lbl[]{3,20;$[\widetilde{c_-}]$}
    \lbl[]{132,39;$[\widetilde{c_2}]$}
    \lbl[]{74,53;$\widetilde{p}$}
    \lbl[]{42,70;$g(\gamma)$}
    \lbl[]{90,75;$\gamma$}
    \lbl[]{84,131;$[\widetilde{c_1}]$}
\end{lpic}
\end{minipage}\hfill
\begin{minipage}{.45\linewidth}
\centering
\vspace{.2cm}
\begin{lpic}{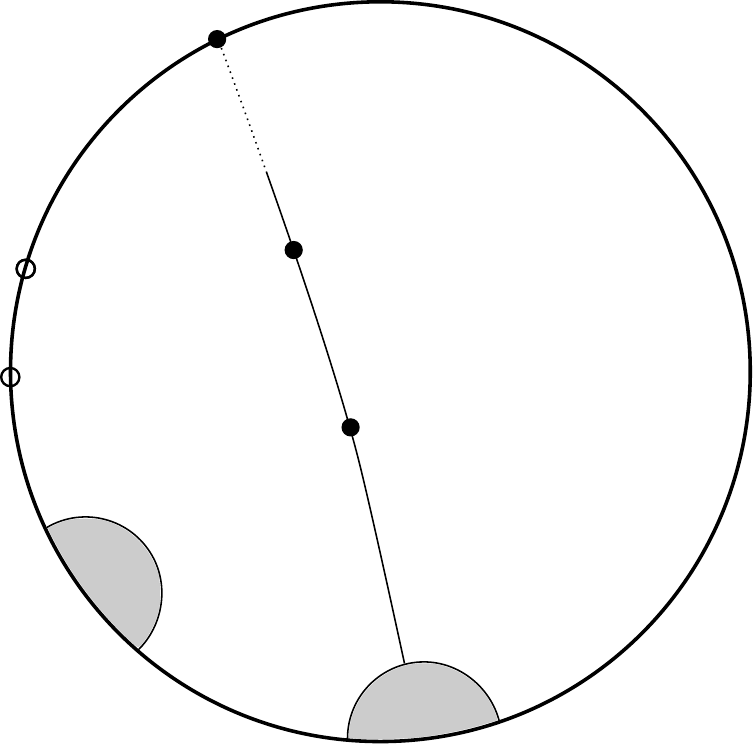(7cm)}
\Large
    \lbl[]{5,22;$[\widetilde{c_-}]$}
    \lbl[]{72,-8;$[g(\widetilde{c_-})]$}
    \lbl[]{71,55;$O$}
    \lbl[]{63,85;$g(O)$}
    \lbl[]{-5,62;$g^-$}
    \lbl[]{-2,82;$g^+$}
    \lbl[]{22,120;$q$}
\end{lpic}
\end{minipage}

\caption{In the translation of \cite[Lem.~2.5]{Przytycki}, tips become intervals on $\partial_\infty\bH$.}
\label{fig:Przytycki}
\end{figure}

Suppose that $O \in g(\gamma)\cap \gamma$, and that without loss of generality $g(O)$ is farther along $g(\gamma)$ from $g(\widetilde{c_-})$ than $O$ (otherwise we can replace $\gamma$ with $g(\gamma)$ and $g$ with $g^{-1}$), and let $q\in \partial_\infty \bH$ be the endpoint of the geodesic extension of $g(\gamma)$. 
Let $H\subset \bH$ be the union of the ray from $g(\widetilde{c_-})$ to $q$ with the complementary halfspace of $g(\widetilde{c_-})$ not containing $q$.

Suppose that $g$ is loxodromic, with fixed points $g^\pm$.
Because $O$ and $g(O)$ are both on $g(\gamma)$, both $g^\pm$ lie to one side of the geodesic containing $g(\gamma)$.
Because $O$ is before $g(O)$ on $g(\gamma)$ and because $[\widetilde{c_-}]$ and $[g(\widetilde{c_-})]$ are disjoint intervals, $[\widetilde{c_-}]$ separates the repelling point $g^-$ from $[g(\widetilde{c_-})]$, and $g^-$ separates $g^+$ from $[\widetilde{c_-}]$, on the component of $\partial_\infty\bH\setminus H$ containing $[\widetilde{c_-}]$. 
Assume without loss of generality that the disjoint intervals $[\widetilde{c_-}]$, $[g(\widetilde{c_-})]$, $q$ are positively oriented on $\partial_\infty \bH$, so that both fixed points $g^\pm$ of $g$ lie in the interval from $q$ to $[\widetilde{c_-}]$. Observe that $g$ now restricts on the interval $[[g(\widetilde{c_-})]^+,g^+]$ to an increasing function.
See Figure~\ref{fig:Przytycki}.

\begin{figure}[ht]
\centering
\begin{lpic}{pics/Przytycki3(7cm)}
\Large
\lbl[]{3,20;$[\widetilde{c_-}]$}
\lbl[]{71,55;$O$}
\lbl[]{68,115;$\widetilde{c_1}$}
\lbl[]{50,23;$\widetilde{c_2}$}
\end{lpic}
\caption{If the tops $\widetilde{c_1}$ and $\widetilde{c_2}$ of $N$ have endpoints in the complementary component of $H$ containing $\widetilde{c_-}$,  stems of $N$ lie to the left of the dotted vertical geodesic.}
\label{fig:Przytycki3}
\end{figure}

If one of the tops $\widetilde{c_i}$ has both endpoints in $[[g(\widetilde{c_-})]^+,g^+]$, then $([\widetilde{c_-}],[\widetilde{c_i}])$ and $([g(\widetilde{c_-})],[g(\widetilde{c_i})])$ are linked pairs, contradicting Lemma~\ref{lem:main sides}.
Therefore we assume that both tops have an endpoint on the component of $\partial_\infty\bH\setminus H$ containing $[\widetilde{c_-}]$.
But then the top arc of the inner flower of $N$ is a geodesic contained in the side of $\bH\setminus H$ containing $\widetilde{c_-}$, and it follows that $\gamma$ and $g(\gamma)$ are disjoint; see Figure~\ref{fig:Przytycki3}.

If instead $g$ is parabolic, the proof above works verbatim, with $g^-=g^+$. \qedhere

\end{proof}

\section{Maximal and type-maximal tulips}
\label{sec:maximal}

\begin{figure}[ht] 
\centering
\includegraphics[width=13cm]{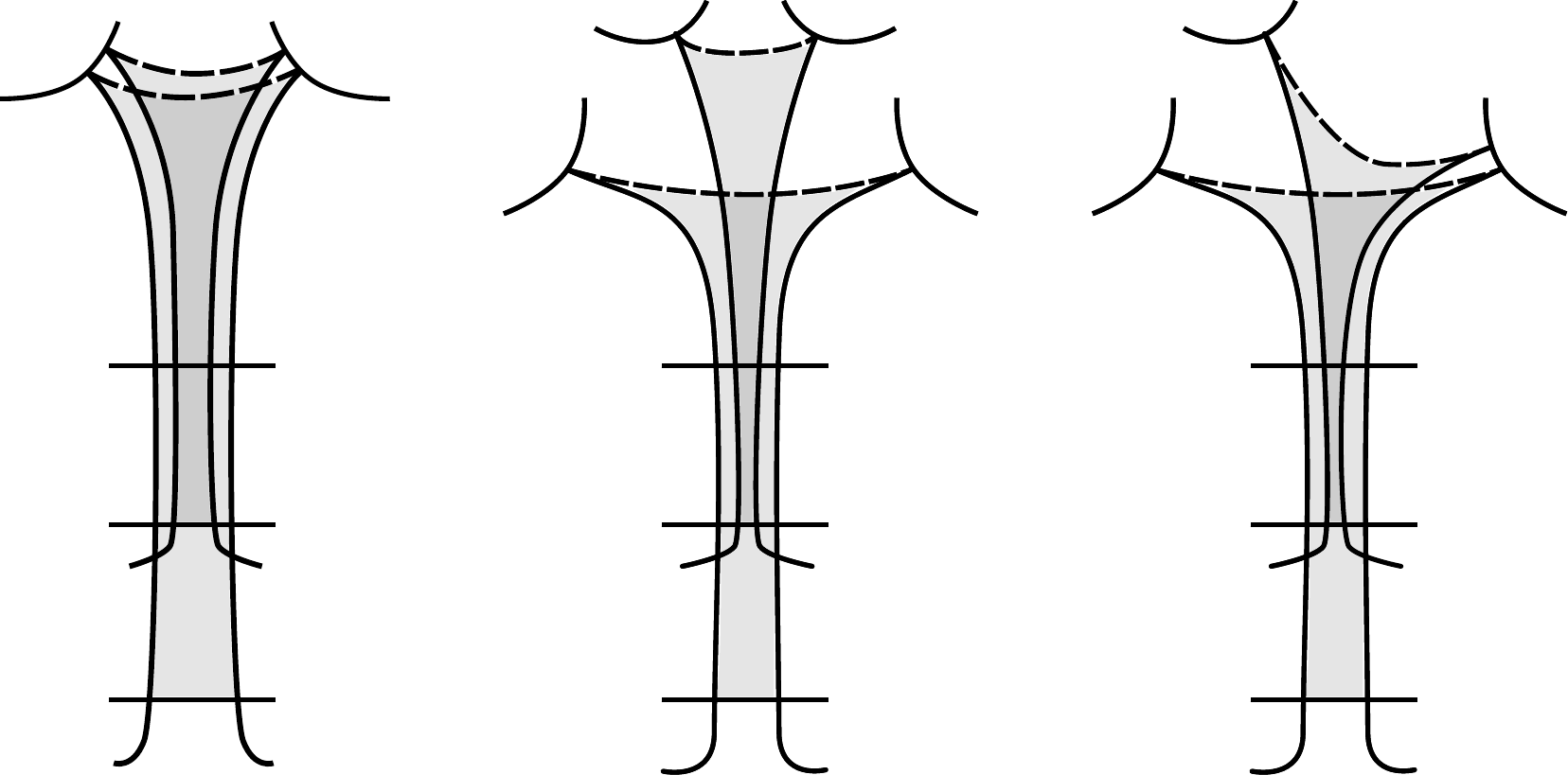}
\caption{Several ways one stem might be a prefix of another.}
\label{Fig:nonmaximal}
\end{figure}

Observe that it is possible for two tulips $N$  and $M$ to have a certain `nestedness': there might be a point $p$ in both inner flowers, so that the stem at $p$ for $N$ is homotopic \emph{rel $\cC$} to a prefix for the stem at $p$ for $M$, see Figure~\ref{Fig:nonmaximal}. It will be necessary to 
avoid this behavior. For simplicity, we will say `$\alpha$ is a prefix of $\beta$' as shorthand for `$\alpha$ is homotopic \emph{rel $\cC$} to a prefix of $\beta$'.

\begin{remark} 
\label{rem:homotopy}
We caution the reader that we use \emph{rel} 
as follows: `$\alpha$ is homotopic to 
$\beta$ \emph{rel $\cC$}' means that 
$\alpha$ and $\beta$ are homotopic \emph{as maps of pairs} $(I,\partial I)\to (S,\{p\}\cup \cC)$.
In particular, we allow one to slide the endpoint of $\alpha$ along $\cC$ during the homotopy.   
\end{remark}

\begin{definition}
Let $N$ and $M$ be distinct tulips, and let $p$ be a point in both inner flowers.
We write $N\prec_p M$ if the $N$-stem at $p$ is a prefix of the $M$-stem at $p$, and $N\prec M$ if there exists a $p$ so that $N\prec_p M$.
A tulip $N$ is \emph{maximal} if it is maximal with respect to $\prec$, i.e.~provided that no stem of $N$ 
at a point in its inner flower is a prefix of a stem in the inner flower of any other tulip $M$.

We write $M\preceq N$ for $M\prec N$ or $M=N$. 
\end{definition}

\begin{remark}
\label{rem:maximal}
We make the trivial observation that the meaning of `maximal' depends on the context of the collection of tulips $\cN$. That is, a tulip $N$ might be maximal in a subset $\cM\subset\cN$, but not maximal in $\cN$.
Note as well that it is tempting to interpret `maximality' for tulips with respect to inclusion. However, the reader is warned that this is not literally true: it is possible for a stem of a tulip $N$ to be the prefix of the stem for tulip $M$, and $N\nsubseteq M$, see Figure~\ref{Fig:nonmaximal}.
\end{remark}

\begin{remark} The reader may also be wondering about the utility of this definition; in particular, one could ask what would go wrong if instead we worked with \emph{all} tulips. As described in 
\S\ref{sec:narration}, our strategy proceeds by bounding the size of a fiber of the projection map from the set of tulips to the surface. So-called \emph{stem systems} (see the beginning of \S\ref{sec:stems} for the formal definition) serve as a stand-in for measuring this fiber. The crux of the issue is that if one does not require a maximality condition as above, it is easy to construct an example of a stem system whose size is \emph{quadratic} in $\chi$ (see Remark \ref{rem: needformaximal}), as opposed to the desired linear-in-$\chi$ bound expressed in Lemma \ref{prop:stembounds}, upon which the main theorem heavily relies. 
\end{remark}

\begin{remark}
Again, the reader is reassured that, though it is not pictured in Figure~\ref{Fig:nonmaximal}, the main curves of a tulip are allowed to intersect, as in Figure~\ref{fig:tulip}.
\end{remark}

Given a tulip $N$ (maximal or otherwise) it will be of use to keep track of the various ways that its main sides behave near 
the tops of $N$. 
Indeed, observe in Figure~\ref{Fig:nonmaximal} 
a certain property a non-maximal tulip might enjoy: in the order at the base, there are outside neighbors of its main curves which \emph{track} its main sides, in that they determine identical trajectories in the dual graph $G(\cP_\cC)$. On the left, these curves track the main sides all the way to the tops of the tulip, while for the two rightmost pictures the tracking fails at the very top. 
It turns out that one of these-- the left one-- is easier to exploit, so we are led to slightly refine our notion of tulip, adding the notion of \emph{type}.

We first observe an obvious consequence of non-maximality: 
when one stem is a prefix of another, the main sides (and main curves) of the larger tulip must track the stalk of the shorter one.

\begin{lemma}
\label{lem:sufficient condition}
Suppose that $N\prec M$. Then the main curves of $M$ track the stalk of $N$, and, 
in the linear order along the lift of the base of $N$, the lifts of main curves of $M$ lie on either side of the lifts of main curves of $N$.
Moreover, any lift of a curve from $\Gamma$ lying between the main curves of $M$, in the order at the base of $N$, must track the stalk of $N$ as well.
\end{lemma}

\begin{proof}
The interiors of the main sides of a tulip intersect the same lifts in $\pi^{-1}\cP_\cC$ as a stem from its inner flower. Hence the main curves of $M$ intersect the same lifts as the stalk of $N$. 
If both tops of $M$ lie to one side of those of $N$, then the inner flowers would be disjoint. 

Indeed, we explain this last sentence for an interested reader: suppose the main curves of $N$ are $\alpha_1$ and $\alpha_2$, with base $c$ and relevant lifts $\widetilde{\alpha_1}$, $\widetilde{\alpha_2}$, and $\widetilde{c}$, so that $\widetilde{\alpha_1},\widetilde{\alpha_2}$ are consecutive in the linear order $\le_{\widetilde{c}}$ on 
$\Gamma_{\widetilde{c}}$. 
Define similarly for $M$ the curves $\beta_1$, $\beta_2$, $d$, and their lifts, and the order $\le_{\widetilde{d}}$ induced by $\widetilde{d}$. Observe that we now know that $\widetilde{\beta_1}$ and $\widetilde{\beta_2}$ both intersect $\widetilde{c}$. Moreover, because $\widetilde{\beta_1}\le_{\widetilde{d}}\widetilde{\beta_2}$, we have $\widetilde{\beta_1}\le_{\widetilde{c}} \widetilde{\beta_2}$ as well. Now there are three possibilities: $\widetilde{\alpha_1}\le_{\widetilde{c}} \widetilde{\alpha_2}\le_{\widetilde{c}}\widetilde{\beta_1}\le_{\widetilde{c}} \widetilde{\beta_2}$, $\widetilde{\beta_1}\le_{\widetilde{c}} \widetilde{\alpha_1}\le_{\widetilde{c}} \widetilde{\alpha_2}\le_{\widetilde{c}}\widetilde{\beta_2}$, or $\widetilde{\beta_1}\le_{\widetilde{c}}\widetilde{\beta_2}\le_{\widetilde{c}}\widetilde{\alpha_1}\le_{\widetilde{c}} \widetilde{\alpha_2}$. The first and third possibilities would imply that the inner flowers are disjoint, a contradiction. Therefore $\widetilde{\beta_1}\le_{\widetilde{c}} \widetilde{\alpha_1}\le_{\widetilde{c}} \widetilde{\alpha_2}\le_{\widetilde{c}}\widetilde{\beta_2}$, as claimed.
\end{proof}

\begin{remark}
Observe in the proof above that it is possible that either one (but not both) of the main curves of $M$ can be a main curve of $N$ as well; for instance, $\widetilde{\beta_1}=\widetilde{\alpha_1}<_{\widetilde{c}}\widetilde{\alpha_2}<_{\widetilde{c}}\widetilde{\beta_2}$. That is, `on either side' in the conclusion of Lemma~\ref{lem:sufficient condition} should be interpreted with `or equal to'.
\end{remark}

\subsection{The type of a tulip}
\label{subsec:type}
We now introduce the notion of \emph{type} for a tulip. The purpose of the type is to make a stronger version of Lemma~\ref{lem:sufficient condition} available (see Lemma~\ref{lem:maximal tops}).

Before diving into the definition, we give an intuitive idea of type: 
the flower of a tulip looks roughly like an ideal triangle, 
mapping to $S$ by a local isometry. 
The hyperbolic surface $S$ has the appearance of almost three-punctured spheres, pieced together along thin annuli around the curves in the pants decomposition $\cP_\cC$. 
The flower maps to a complementary pants in one of only a handful of ways, corresponding to how 
the main curves spin around the base of the flower, a 
cuff of the pants,
since a pair of curves cannot twist an arbitrarily different number of times without intersecting more than once.
The type indexes possible configurations for the flower, see Figure~\ref{fig:types}.

\begin{figure}
\vspace{1cm}
\begin{lpic}{pics/nibTypesAgain(12cm)}
\LARGE
\lbl[]{400,240;$d=1$}
\lbl[]{425,60;$d=\frac12$}
\end{lpic}
\caption{Two of the seven possible types of flowering for tulips.}
\label{fig:types}
\end{figure}

The key point is that when a pair of tulips are of the same type and a stem of one is a prefix of the other
, the two configurations on the right of Figure~\ref{Fig:nonmaximal} cannot occur -- this is the content of Lemma \ref{lem:maximal tops} below. 
One gains an extra bit of control: 
when it comes to a pair of tulips of the same type the ``tracking'' behavior cannot fail at the very top. 
This enables one to learn something by `looking next to a tulip' at its base, the crucial ingredient in the proof of Lemma~\ref{lem:accounting}, which is essentially the entire proof of Theorem~\ref{thm:P2}.

Now, let $\alpha, \omega$ denote the main sides of a tulip $N$, let $c \in \mathcal{P}_{\cC}$ be the 
base of the flower of $N$, and let $P$ be the pair of pants, 
complementary to $\cP_\cC$, containing the flower of $N$.
Evidently one component of $\partial P$ is  $c$; suppose the other two are  $c_1,c_2\in\cP_\cC$.

Fix lifts $\tilde{c}, \widetilde{\alpha}, \widetilde{\omega}$ of $c, \alpha, \omega$ to the universal cover as described in the construction of tulips; associated to them is a fixed lift $\widetilde{P}$ of the universal cover of $P$ and a vertex $v$ in the dual graph $G(\cP_\cC)$ to $\pi^{-1} \cP_\cC$. 
Observe that simple arcs on a pair of pants are quite restricted. For instance, there are only two homotopy classes of simple arcs with one endpoint on $c$ and one endpoint on $c_1$ or $c_2$.
It follows that any lift of a simple closed curve intersecting $c$ once determines a ray in $G(\cP_\cC)$, in such a way that the edge immediately following $v$ comes from a discrete set, naturally parameterized by 
\[
A_0=\{\ldots,\; g^{-1}\cdot\widetilde{c_1}\;,\; g^{-1}\cdot\widetilde{c_2}\;,\;\widetilde{c_1}\;,\;\widetilde{c_2}\;,\;g\cdot\widetilde{c_1}\;,\;g\cdot\widetilde{c_2}\;,\ldots\}\approx \bZ~,
\]
where $g\in\pi_1S$ is a primitive covering transformation with translation axis $\tilde{c}$, and $\widetilde{c_1}$ and $\widetilde{c_2}$ are lifts of $c_1$ and $c_2$ on the boundary of $\widetilde{P}$, incident to $v$. 

Slightly more generally, there are \emph{four} homotopy classes \emph{rel $\cP_\cC$} of simple (oriented) arcs on $P$ with one endpoint on $c$, since the simple arc might have both of its endpoints on $c$. The lifts of such arcs determine trajectories in $G(\cP_\cC)$ where the edge following $v$ comes from the set
\begin{equation}
\label{eq:A}
A=A_0 \cup 
\{\ldots,\; g^{-1}\cdot\widetilde{c_0}\;,\; g^{-1}\cdot\widetilde{c_0}'\;,\;\widetilde{c_0}\;,\;\widetilde{c_0}'\;,\;g\cdot\widetilde{c_0}\;,\;g\cdot\widetilde{c_0}'\;,\ldots\} 
~,
\end{equation}
where $\widetilde{c_0}$ and $\widetilde{c_0}'$ are lifts of $c$ on $\partial\widetilde{P}$. 
Observe that components of $\partial\widetilde{P}$ are linearly ordered from the  the viewpoint of $\tilde{c}$,
and 
choices for lifts can be made so that $\widetilde{c_1},\widetilde{c_0},\widetilde{c_2},\widetilde{c_0}'$ are consecutive in the linear order along $\tilde{c}$.
One finds an order-preserving bijection $A\approx \bZ[1/2]$, well-defined up to isometries of $\bZ[1/2]$.
We define the 
\emph{annular distance} $d_A(x,y)$ between $x,y\in A$ to be $|x-y|$ (with the identification $A\approx \bZ[1/2]$ implicit).

Now $\tilde{\alpha}$ and $\tilde{\omega}$ are lifts of simple curves intersecting $\tilde{c}$, so they each determine edges $e_\alpha,e_\omega\in A$.

\begin{definition} \label{def:type} The \textit{type} of $N$ is the annular  distance 
$d_A(e_\alpha,e_\omega)$. 
\end{definition}

\begin{lemma}
\label{lem:enumerate types}
The type of a tulip satisfies $d_A(e_\alpha,e_\omega)\in \{\frac12,1,\frac32,2,\frac52,3,4\}$.
\end{lemma}

\begin{proof}
Observe that $d_A(e_\alpha,e_\omega)>0$, because by definition the main curves $\widetilde{\alpha}$ and $\widetilde{\omega}$ have endpoints on distinct components of $\partial\widetilde{P}$.
Moreover, it is not hard to see that $d_A(e_\alpha,e_\omega)\le 4$: otherwise the images of $\widetilde{\alpha}$ and $\widetilde{\omega}$ would intersect at least twice.
Finally, we observe that $d_A(e_\alpha,e_\omega)=\frac72$ is not possible. Indeed, up to covering transformations we would find endpoints $\{0,\frac72\}$, $\{1,\frac92\}$, $\{\frac12,4\}$, or $\{\frac32,5\}$, and it is straightforward to check that the main curves would intersect at least twice.
\end{proof}

\begin{lemma}
\label{lem:maximal tops}
Suppose that 
$N$ and $M$ are of the same type, and that $N\prec M$. Then the main sides of $M$ track the main sides of $N$, from the base of $N$ to its tops.
\end{lemma}

\begin{proof}
We adopt notation for $N$: the main curves $\alpha_1$, $\alpha_2$ have lifts $\widetilde{\alpha_1}$, $\widetilde{\alpha_2}$ intersecting $\widetilde{c_-}$, a lift of $c_-$, the base of $N$; 
the tops of $N$ lie on $\widetilde{c_1},\widetilde{c_2}\in A$.
Because 
$N\prec M$, by Lemma~\ref{lem:sufficient condition} the main curves $\widetilde{\beta_1}$,$\widetilde{\beta_2}$ of $M$ track the stalk of $N$.
In particular, the base of the flower for $M$ is 
the same as that of $N$; suppose that the tops of $M$ occur at $\widetilde{d_1},\widetilde{d_2}\in A$. 

Now consider the relative positions of $\widetilde{\alpha_1}$,  $\widetilde{\alpha_2}$,  $\widetilde{\beta_1}$, $\widetilde{\beta_2}$ in the order at $\widetilde{c_-}$. 
By Lemma~\ref{lem:sufficient condition}, 
one top of $M$ is to the left of $\widetilde{c_1}$ and the other is to the right of $\widetilde{c_2}$.
Because $N$ and $M$ are of the same type, $d_A(\widetilde{c_1},\widetilde{c_2})=d_A(\widetilde{d_1},\widetilde{d_2})$, so we must have $\{\widetilde{c_1},\widetilde{c_2}\}=\{\widetilde{d_1},\widetilde{d_2}\}$. This means that the main sides of $M$ track the main sides of $N$ from the base of $N$
to its tops, as claimed.
\end{proof}

With our new notion of 
type in hand, we consider a type-specific version of maximality: 

\begin{definition} A tulip of 
type $d$ is \emph{type-maximal} if 
it is maximal amongst the tulips of type $d$.
\end{definition}

\begin{lemma}
\label{lem:maximal exists}
Every tulip of type $d$ is $\preceq$-dominated by a type-maximal tulip of type $d$.
\end{lemma}

\begin{proof}
Recall that $\epsilon$ has been chosen small enough that $N$ has nonempty inner flower. If $N$ is not type-maximal, then $N\prec M$ for some tulip $M$ of the same type. In that case, Lemma~\ref{lem:maximal tops} implies that $N$ and $M$ have identical inner flowers. 
We conclude that $\prec$ is transitive on tulips of fixed type.

Finally, we claim that no $\prec$-chain can have length more than $\max\{ \iota(\gamma,\cP_\cC):\gamma\in \Gamma\}$: observe that $N\prec M$ implies that the combinatorial length of the stalk of $N$ is strictly less than that of $M$, but each stalk may intersect at most $\max\{ \iota(\gamma,\cP_\cC):\gamma\in \Gamma\}$ pants lifts. 
\end{proof}

\begin{remark}
It follows that $\prec$ is a strict partial order on tulips of a given type. 
In fact, $\prec$ may not be a strict partial order on all of the tulips, unrestricted by type, precisely because transitivity of $\prec$ may fail: it is possible that $N_1$, $N_2$ and $N_2$, $N_3$ have overlapping inner flowers, but those of $N_1$ and $N_3$ are disjoint. Hence it is possible that $N_1\prec N_2$ and $N_2\prec N_3$, but $N_1\nprec N_3$. 
\end{remark}

\subsection{Curve accounting with type-maximal tulips}
Here we demonstrate that every curve $\gamma \in \Gamma$ is the main side of a type-maximal tulip. 
Later, we show that all type-maximal tulips of fixed type has a cardinality that is bounded above by an explicit quadratic function of $\chi$, and the proof of the main theorem will then follow from Lemma~\ref{lem:enumerate types}.

\begin{lemma} \label{lem:accounting} Let $\gamma \in \Gamma\setminus\cC$, and let $\Vec{\gamma}$ be an oriented representative for $\gamma$. Then $\Vec{\gamma}$ is a main curve of a tulip that is type-maximal. 
\end{lemma}

\begin{proof} 
Given $\gamma$ as above (but suppressing the 
orientation notation), let $\omega \in \Gamma$ be some other curve that spans a tulip $N$ with $\gamma$. 
By Lemma~\ref{lem:maximal exists} there is a type-maximal tulip $M$ 
so that $N\preceq M$.
Let $c \in \cC$ be the curve containing the base of $N$, and let $\widetilde{\gamma}, \widetilde{c}, \widetilde{\omega}$ be the relevant lifts for $N$. 
Let $\le$ indicate the linear order induced by $\widetilde{c}$ on $\Gamma_{\widetilde{c}}$.
We may assume without loss of generality that $\widetilde{\gamma}< \widetilde{\omega}$.

\begin{figure}
    \centering
\begin{lpic}{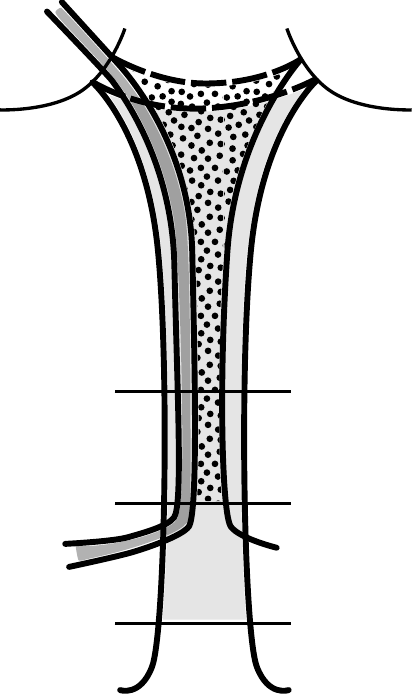(,9cm)}
\Large
\lbl[]{-5,100;$\widetilde{x}$}
\lbl[]{4,28;$\widetilde{\zeta}$}
\lbl[]{6,17;$\widetilde{\gamma}$}
\lbl[]{13,2;$\widetilde{\delta}$}
\lbl[]{51,22;$\widetilde{\omega}$}
\lbl[]{53,33;$\widetilde{c}$}
\end{lpic}
\caption{A schematic of the proof of Lemma~\ref{lem:accounting}: The non-type-maximal tulip $N$ is dotted; $M$, lightly shaded, is type-maximal and $N\prec M$; $M_1$, only partially shown, is darkly shaded; $M_2$, not pictured, is type-maximal with $M_1\preceq M_2$. 
Because no curve lies strictly between $\widetilde{\gamma}$ and $\widetilde{\omega}$, $\widetilde{\gamma}$ must be a main curve of $M_2$.
}
\label{fig:nonmaximalAdjacent}
\end{figure}

If $\widetilde{\gamma}$ is one of the main curves of $M$, there is 
nothing to prove. Otherwise,
by Lemma~\ref{lem:sufficient condition} there must be a main curve $
\widetilde{\delta}$ of $M$ so that $\widetilde{\delta}<\widetilde{\gamma}$.
Moreover, Lemma~\ref{lem:maximal tops} implies that $\widetilde{\delta}$ terminates at the same lift $\widetilde{x}\in A$ that $\widetilde{\gamma}$ does, and that it fellow-travels with $\widetilde{\gamma}$ 
from $\widetilde{c}$ to $\widetilde{x}$, in that they induce the same path in $G(\cP_{\cC})$. 
By Lemma~\ref{lem:discrete}, there is an immediate predecessor $\widetilde{\zeta}$ of $\widetilde{\gamma}$ along $\widetilde{c}$.
Hence $\widetilde{\gamma}$ and $\widetilde{\zeta}$ are the main curves of some tulip $M_1$.
Now 
we have $\widetilde{\delta}\le \widetilde{\zeta}<\widetilde{\gamma}$, 
and, because the geodesics $\widetilde{\gamma}$ and $\widetilde{\delta}$ fellow-travel from $\widetilde{c}$ to $\widetilde{x}$ in $G(\cP_\cC)$,
so $\widetilde{\zeta}$ does as well. 
See Figure~\ref{fig:nonmaximalAdjacent}.

Again by Lemma~\ref{lem:maximal exists}, there exists a type-maximal tulip $M_2$ so that $M_1\preceq M_2$. We claim that $\widetilde{\gamma}$ is a main curve of $M_2$, completing the proof.
Indeed, if $\widetilde{\gamma}$ is not a main curve for $M_2$, then by Lemma~\ref{lem:sufficient condition} a main curve of $M_2$, say $\widetilde{\eta}$, lies to the right of $\widetilde{\gamma}$ at the base of $M_2$ 
and fellow-travels the stalk of $M_1$, which means that $\widetilde{\eta}$ fellow-travels $\widetilde{\gamma}$ from the base of $M_1$ to $\widetilde{x}$ (in fact, to the base of the flower of $M_1$). 
This means that 
$\widetilde{\gamma}<\widetilde{\eta}<\widetilde{\omega}$, a contradiction.

To summarize: either the type-maximal tulip $M$ that $\preceq$-dominates $N$ has $\widetilde{\gamma}$ as a main curve, or else there is a tulip $M_1$, immediately on the other side (along $\widetilde{c}$) of $N$ from $\widetilde{\gamma}$, and the type-maximal tulip $M_2$ that $\preceq$-dominates $M_1$ must have $\widetilde{\gamma}$ as a main curve.
\end{proof}

\begin{remark}
Again, we reassure the reader that the main curves of a tulip may cross as in Figure~\ref{fig:tulip}, though this is not pictured in Figure~\ref{fig:nonmaximalAdjacent}. The proof of Lemma~\ref{lem:accounting}, like the other proofs in this section, is agnostic to this distinction.
\end{remark}

\begin{remark} \label{type can change} 
%
At the last step of the proof, Lemma~\ref{lem:sufficient condition} implies that $M_2$ is actually a \emph{maximal} tulip, not merely type-maximal. This makes it seem like the conclusion of the lemma can be promoted to a stronger one. However, 
this proof required at the beginning that the tulip $N$ was \emph{not type-maximal}. 
So there is no apparent contradiction if $\gamma$ is the side of a tulip that is only type-maximal. 
Moreover, the tulips $M_1$ and $M_2$ 
above may be of different types than the tulip $N$ (or than each other), so we have no control over the type of the type-maximal tulip that accounts for $\gamma$.

We are unable to determine whether every curve $\gamma\in \Gamma$ is the main curve of a maximal tulip.
\end{remark}

\subsection{The refined type of a tulip}
\label{subsec:restricted type}
Before we turn to analyzing fibers of the map from type-maximal tulips to the surface, we include one more refinement that will simplify the coming proofs.

For each pair of pants component $P$ of $S\setminus \cP_\cC$, 
fix 
a homeomorphism of its completion $\overline{P_0}$ with $P_0$, 
where $P_0\subset S^2$ is the complement of three disks $D_1$, $D_2$, $D_3$, with disjoint closures.

For each tulip $N\to S$, the image of its flower is contained in a pair of pants $P$, and the base of its flower is one of the cuffs of $P$. 
With the identification $P\approx P_0$ implicit, the base of the flower of $N$ is incident to $\partial D_i$, for some $i\in\{1,2,3\}$. 


\begin{definition}
\label{def:refined type}
The \emph{refined type} of a tulip $N\to S$ is the pair $(d,k)$, where $N$ is of type $d$, and 
the base of the flower of $N$ is incident to $D_k$.
\end{definition}

Observe that Lemma~\ref{lem:enumerate types} has the following immediate consequence:

\begin{lemma}
\label{lem:refined types}
There are at most 21 refined types of tulips associated to a 1-system.
\end{lemma}

Going forward, we will say that the flower of a tulip \emph{blooms over $(c,P)$}, with $c \in \cC$ and $P$ a pair of pants in the complement of $\cC$, when the base of the flower is a lift of $c$ and when the projection of the flower lives in $P$; when $P$ is clear from context, we will say that the floor blooms over $c$.

By examination of flowers of each of the seven types given in Lemma~\ref{lem:enumerate types}, we can see that each inner flower arising from a 1-system maps to a complementary pants with fibers of size at most two. It follows that, for each refined type and for each point in one of the complementary pants, there are only two possible homotopy classes of arcs to the specified cuff that a stem of that type begins with (recall that the endpoint of the arc is allowed to move along the cuff during the homotopy).
In fact, one can observe that these two homotopy classes are independent of the type. Therefore:

\begin{lemma}
\label{lem:homotopy classes}
Let $P$ be a pair of pants complementary to $\cP_\cC$, let $p\in P$, and let $c$ be a cuff of $P$.
There are two arcs $\alpha_1$ and $\alpha_2$, from $p$ to $c$, with the property:
For any tulip $N\to S$ of a refined type so that the flower blooms over $(c,P)$, and so that $p$ lies in the inner flower with stem $\alpha$ at $p$, there is a prefix $\alpha_0\subset \alpha$ so that $\alpha_0 \simeq \alpha_i$ rel $\cC$, for $i\in\{1,2\}$. 
\end{lemma}

We note that we will only need the weaker statement: 
for each refined type and each point in a pair of pants, there are a pair of homotopy classes of arcs containing a prefix of every stem at the point of a tulip whose flower blooms over the given pants cuff. 

We are now ready to proceed by bounding the size of the set of all 
tulips that are maximal amongst their type.

\section{Bounds for stem systems} \label{sec:stems}

Let $S$ be a hyperbolic surface, let $\cC$ be a geodesic multicurve, contained in a pants decomposition all of whose cuffs have length $\epsilon$ (a small quantity, soon to be sent to zero). 
The main result of this section is that the cardinality of a stem system is bounded by a linear function of $|\chi|$.

\begin{definition}[Stem system on $S$]
\label{def:stem system}
Let $S$ and $\cC$ be as above, and let $p\in S\setminus\cC$. A \emph{stem system} based at $(p,\cC)$ is a collection $\cA$ of simple geodesic arcs from $p$ to $\cC$ satisfying, for all $\alpha,\beta\in\cA$:
\begin{enumerate}[label=(\roman*)]
\item \label{it:stem prefix}
If there is a prefix $\beta_0\subset\beta$ homotopic to $\alpha$ \emph{rel $\cC$}, then $\alpha=\beta$.
\item 
\label{it:stem triangles}
If $\alpha$ and $\beta$ intersect, they do so only once, transversely, and the intersection point is on the boundary of a triangle formed by sub-arcs of $\alpha$, $\beta$, and one of their base curves in $\cC$. 
\end{enumerate}
\end{definition}
Recall that `$\alpha$ is homotopic \emph{rel $\cC$} to $\beta$' means that the points in $\cC$ are allowed to slide during the homotopy, see Remark~\ref{rem:homotopy}.

We give a rough sketch of the proof: If there are many stems emanating from a point $p$, then the counterclockwise order around $p$ arranges them into triangular pieces, each of which maps to the surface in a manner reminiscent of the tulips from the original $1$-system--these are something like \emph{secondary tulips}. 
Just as in \S\ref{sec:maximal},
as $\epsilon$ approaches zero the geometry of each complementary pants simplifies,
and one finds that the flowers of these secondary tulips map to the surface in one of a handful of highly constrained ways--these are the types.
There is an important bit of added benefit: there is a uniform bound on the number of inner flowers of secondary tulips that are of the same type, blooming at the same cuff; see Lemma~\ref{lem:flowers in same position}.
This leads to the needed linear-in-$|\chi|$ bound.

\begin{proposition} \label{prop:stembounds} 
If $\cA$ is a stem system based at $(p,\cC)$, then $|\cA|\leq 24|\chi|+4$.
\end{proposition}

\begin{proof}
Observe that neither condition of Definition~\ref{def:stem system} is affected by adding disjoint curves to $\cC$, so we may safely assume $\cC$ is a pants decomposition below.
Choose a lift $\widetilde{p}\in \bH$, and incident lifts of the stems in $\cA$ to geodesics $\widetilde{\alpha}_1,\ldots,\widetilde{\alpha}_n$ in $\bH$, arranged in cyclic order around $\widetilde{p}$, with endpoints on lifts of curves $\widetilde{c}_1,\ldots,\widetilde{c}_n$ (components of $\pi^{-1}\cC$).
Because of property \ref{it:stem prefix} in the definition above, none of these base curves for the listed stems may be nested from the viewpoint of $\widetilde{p}$. That is, for each $i\ne j$, $\widetilde{p}$ and $\widetilde{c}_i$ lie on the same side of $\widetilde{c}_j$.
Let $\Omega\subset \bH$ indicate the polygon whose vertices are the endpoints $\widetilde{\alpha}_1\cap \widetilde{c}_1,\ldots,\widetilde{\alpha}_n\cap \widetilde{c}_n$. See Figure~\ref{fig:geodesic stem system polygon} for an image in which $\widetilde{p}$ is interior to $\Omega$--the reader may check that this is not necessary in our proof.

By the Collar Lemma, the distances between $\widetilde{c}_i$ and $\widetilde{c}_{i+1}$ are quite large, and it is clear that $\Omega$ approaches an ideal polygon. Roughly speaking, the remainder of the proof demonstrates: 
\[
\emph{Amost all of $\Omega$ almost embeds on $S$ under the covering projection.}
\]
The conclusion is within reach:
$
(n-2)\pi \approx \Area(\Omega)\lesssim \Area(S)=2\pi |\chi|
$.

It remains to make this argument rigorous.
In step 1, we will identify a subset $\Omega_0\subset \Omega$, so that $\Area(\Omega_0)$ is roughly $\Area(\Omega)$ as $\epsilon\to 0$, and in step 2 we show that $\Omega_0$ maps to the surface with uniformly controlled fibers. 
The subset $\Omega_0$ will be a union of triangles, essentially one for each side of $\Omega$ (up to uniform additive error), each limiting geometrically to an ideal triangle as $\epsilon$ goes to $0$, making the computation of $\Area(\Omega_0)$ transparent in the limit. 
By mimicking the type analysis of \S\ref{subsec:type}, we will show that  $\pi|_{\Omega_0}:\Omega_0\to S$ maps at most twenty-four triangles to each pair of pants.

\begin{figure}
\begin{minipage}{.45\textwidth}
\centering
\begin{lpic}{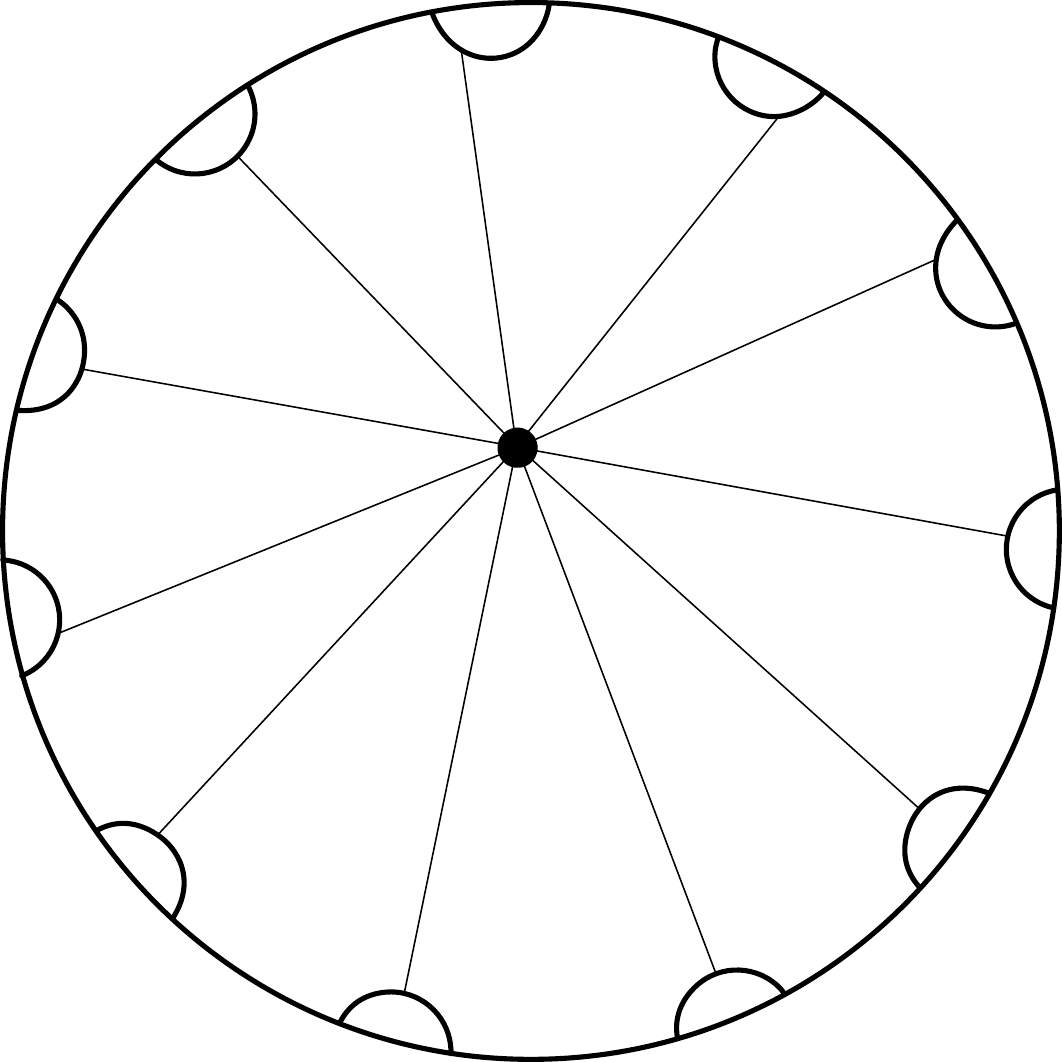(,7cm)}
\Large
\lbl[]{80,120;$\widetilde{p}$}
\lbl[]{140,102;$\widetilde{\alpha}_{i}$}
\lbl[]{164,76;$\widetilde{c}_{i}$}
\lbl[]{130,80;$\widetilde{\alpha}_{i+1}$}
\lbl[]{140,37;$\widetilde{c}_{i+1}$}
\end{lpic}
\end{minipage}\hfill
\begin{minipage}{.45\textwidth}
\centering
\begin{lpic}{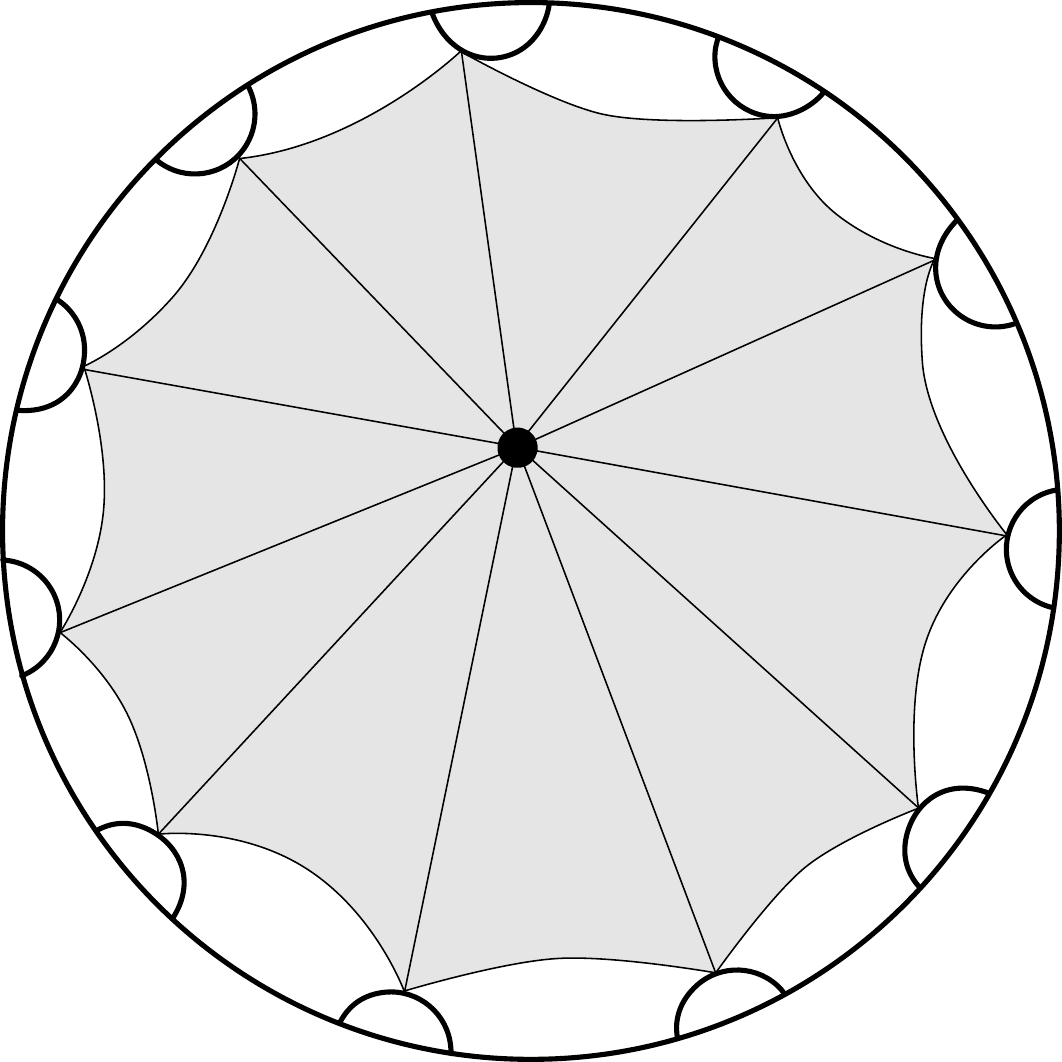(,7cm)}
\Large
\lbl[]{90,40;$\Omega$}
\end{lpic}
\end{minipage}
\caption{A geodesic stem system at $\widetilde{p}$ in the universal cover gives rise to a polygon $\Omega$ which looks more and more ideal as $\epsilon\to 0$.}
\label{fig:geodesic stem system polygon}
\end{figure}

\textbf{Step 1.}  
For each $i$, consider the triangle 
\[
\Delta_i:=\Delta(\widetilde{\alpha}_i\cap \widetilde{c}_i,\widetilde{\alpha}_{i+1}\cap \widetilde{c}_{i+1},\widetilde{p})~.
\]
The intersection of this triangle with $\pi^{-1}\cC$ consists of a (potentially empty) set of geodesic arcs cutting off each of the three corners; at each corner there is an innermost lift from $\pi^{-1}\cC$. 
Let $F_i$ indicate the component of the complement of $\pi^{-1}\cC$ in $\Delta_i$ that is adjacent to the innermost lifts. 
Finally, the triangle $T_i\subset F_i$ is the result of deleting the convex hulls of adjacent innermost lifts from $\pi^{-1}\cC$;
see Figure~\ref{fig:adjust triangles}.
We will refer to $\Delta_i$ as a `secondary tulip', to $\widetilde{\alpha}_i$ and $\widetilde{\alpha}_{i+1}$ as its `stem main sides', to $F_i$ as its `flower', and to $T_i$ as its `inner flower', though the reader is cautioned that these terms are not meant to be literal as in Definition~\ref{def:nibs etc}, but only evocative; see Remark~\ref{rem:secondary nibs}.

Observe that each of $F_i$ has a `base', which is either a sub-arc in some lift in $\pi^{-1}\cC$, or a corner at $\widetilde{p}$.
After relabelling the flowers, we may assume that $F_1,\ldots,F_m$ are \emph{not} $\widetilde{p}$-adjacent.
We claim that $m \ge n-4$. To see this, note that associated to the flower $F_i$ (for $m < i \le n$) is a pair of lifts in $\pi^{-1} \cC$: these are the first curves in $\pi^{-1} \cC$ seen when we traverse either stem side of $F_i$, starting from $\widetilde{p}$. Let us call these lifts $\widetilde{d_{i+1}}, \widetilde{d_i}$; by the assumption that $F_i$ is $\widetilde{p}$-adjacent, $\widetilde{d_{i+1}}$ crosses \textit{only} $\widetilde{\alpha_{i+1}}$ and not $\widetilde{\alpha_{i}}$, and similarly for $\widetilde{d_{i+1}}$. Finally, let $\widetilde{x_{i}}$ be the point of intersection between $\widetilde{\alpha_{i}}$ and $\widetilde{d_{i}}$ and let $\beta_{i} \subset \alpha_{i}$ denote the sub-arc of $\alpha_{i}$ starting at $p$ and ending at $\pi(\widetilde{x_{i}})= x_{i}$.

Each $\widetilde{d_i}$ projects to a cuff of the pair of pants $P$ in the complement of $\cC$ containing $p$. Moreover, each arc $\beta_{i}$ is contained in $P$, and for $i \neq j$, $\beta_i$ can not be homotopic \textit{rel} $\cC$ to $\beta_j$ since otherwise $\widetilde{x_i}$ and $\widetilde{x_{j}}$ would belong to the same lift $\widetilde{d_{i}}$. If the number of such $\widetilde{d_{i}}$'s is at least $5$, we see $5$ embedded arcs $\beta_{m+1},..., \beta_{m+5}$ in $P$, starting at $p$ and ending on cuffs of $P$, no two of which are homotopic \textit{rel} $\cC$. But then at least two of them must intersect at some point other than $p$ in such a way that \ref{it:stem triangles} is violated. Therefore there are at most $4$ such arcs. The stem sides of each flower that is $\widetilde{p}$-adjacent are $\widetilde{d_{j}}, \widetilde{d_{j+1}}$ and such a flower is uniquely determined by these sides. Therefore there are also at most $4$ $\widetilde{p}$-adjacent flowers, as claimed.

\begin{figure}
\begin{minipage}{.45\textwidth}
\centering
\begin{lpic}{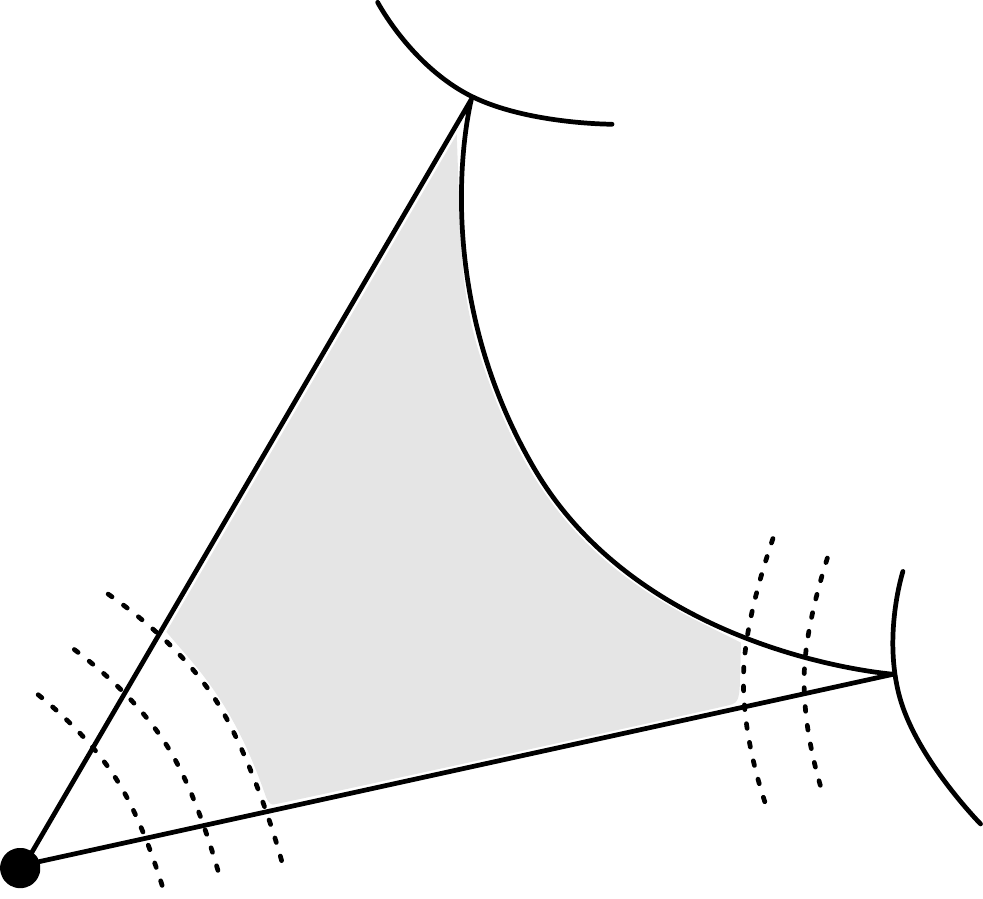(,6.5cm)}
\Large
\lbl[]{-15,5;$\widetilde{p}$}
\lbl[]{40,100;$\widetilde{\alpha}_{i+1}$}
\lbl[]{100,142;$\widetilde{c}_{i+1}$}
\lbl[]{170,35;$\widetilde{c}_i$}
\lbl[]{80,12;$\widetilde{\alpha}_i$}
\lbl[]{70,55;$F_i$}
\end{lpic}
\end{minipage}\hfill
\begin{minipage}{.45\textwidth}
\centering
\begin{lpic}{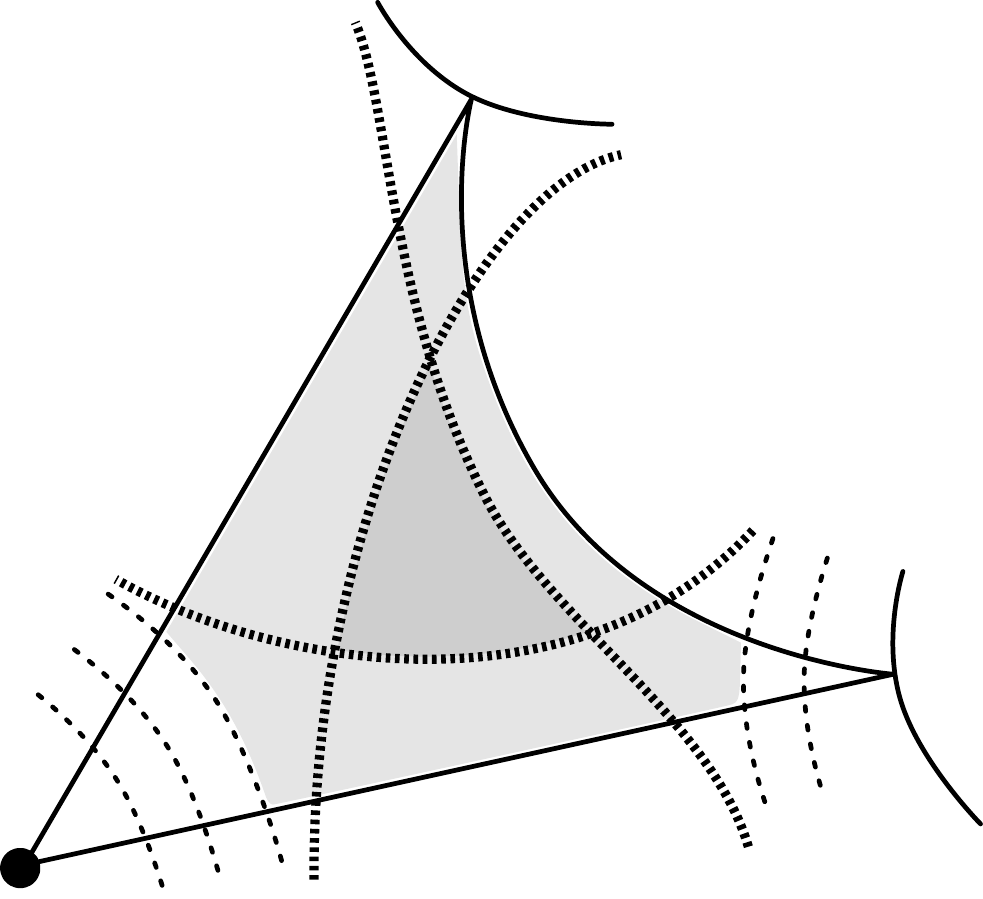(,6.5cm)}
\Large
\lbl[]{76,58;$T_i$}
\end{lpic}
\end{minipage}
\caption{The `flower' $F_i$ is the innermost subset of $\Delta_i$ complementary to $\pi^{-1}\cC$, and $T_i$
is the complement in $F_i$ of convex hulls of adjacent innermost lifts from $\pi^{-1}\cC$.}
\label{fig:adjust triangles}
\end{figure}

Now we define
\[
\Omega_0=\bigsqcup_{i=1}^m T_i~.
\]
As $\epsilon$ goes to $0$, $\Area(\Omega_0)$ approaches $m\pi$. 

\textbf{Step 2.} 
It is possible that $\pi|_{\Omega_0}:\Omega_0\to S$ is not an embedding. Nonetheless, we claim that $\pi|_{\Omega_0}$ sends at most $24$ flowers to each pair of pants.

Consider a pair of pants $P\subset S\setminus \cC$, and suppose that $F_i$ (and, hence, $T_i$) maps to $P$ under $\pi$. 
Observe that this means that the base of $F_i$ occurs at a cuff of $P$, which we denote $c$.

The stem main sides of $F_i$ project to simple arcs on $P$. 
Simple arcs in $P$ are constrained as in \S\ref{subsec:type}. In fact, employing the same notation as in Lemma~\ref{lem:enumerate types}, the annular distance between the stem main sides at $c$ lies in $\{\frac12,1,\frac32,2\}$:
If this distance were larger than two, we would find an intersection 
of stems that does not form a triangle with either base, in violation of \ref{it:stem triangles}. 
As in \S\ref{subsec:type}, we will refer to the annular distance between the stem main sides of the flower $F_i$ as the \emph{type} of $F_i$.

While there are at most four types of inner flowers mapping to each of the three cuffs of each pair of pants, it is possible for more than one flower of a given type to map to the same pair of pants with base at the same cuff.
Nonetheless, more precise control is possible:
\begin{lemma}
\label{lem:flowers in same position}
Suppose that $\cF\subset\{F_i\}$ is a subset of the flowers, each of which maps to $P$ with the same base cuff, and is of type $d$. 
\begin{enumerate}[label=(\roman*)]
    \item If $d=1/2$ then $|\cF|\le 4$.
    \item If $d=1$ then $|\cF|\le 2$.
    \item If $d=3/2$ then $|\cF|\le 1$.
    \item If $d=2$ then $|\cF|\le 1$.
\end{enumerate}
\end{lemma}

Each pair of pants has three cuffs, 
so this lemma implies that at most $3\times(4+2+1+1)=24$ flowers among $\{F_i\}$ can map to each pair of pants.
There are $|\chi|$ pairs of pants on $S$, so we find that $m\le24|\chi|$, whence $n\le 24|\chi|+4.$\qedhere
\end{proof}

Towards Lemma~\ref{lem:flowers in same position}, we make a preliminary observation: if two based inner flowers are in the same $\pi_1S$-orbit, 
the corresponding secondary tulips are in fact the same.

\begin{lemma}
\label{lem:flowers equal}
Suppose that there is $g\in\pi_1S$ so that $g(T_j)=T_i$, and so that $g$ takes the base of $F_j$ to that of $F_i$. Then $g=id$ and $i=j$.
\end{lemma}

\begin{proof}
The assumptions on $g$ imply that $g(\widetilde{c}_j)=\widetilde{c}_i$ and $g(\widetilde{c_{j+1}})=\widetilde{c_{i+1}}$. Indeed, the edge of $T_{j}$ whose geodesic extension has endpoints $[\widetilde{c_{j+1}}]^{-}$ and $[\widetilde{c_{j}}]^{+}$ is sent by $g$ to the edge of $T_{i}$ whose geodesic extension has endpoints $[\widetilde{c_{i+1}}]^{-}$ and $[\widetilde{c_{i}}]^{+}$. Another application of [\cite{Farb-Margalit}, p. 32] implies that these geodesics coincide.

Let $\widetilde{c}$ be the base of $F_i$, and
let $s_i$ and $s_{i+1}$ be the complete geodesics pictured in Figure~\ref{fig:rule out}: both contain $\widetilde{p}$ and pass through the 
endpoints $[\widetilde{c_i}]^-$ and $[\widetilde{c}_{i+1}]^+$ of $\widetilde{c}_i$ and $\widetilde{c}_{i+1}$. 
We will use the definition of the stem system to rule out potential locations for $g(\widetilde{p})$.

We next claim that the assumption that $g$ takes the base of $F_j$ to that of $F_i$ implies that $g(\widetilde{p})$ lies on the same side of $\widetilde{c}$ as $\widetilde{p}$. To see this, note that the base of $F_{j}$ separates $\widetilde{p}$ from $T_{j}$, and thus the base of $g(F_j)= \widetilde{c}$ separates $g(\widetilde{p})$ from $g(T_{j}) = T_i$. Since $\widetilde{c}$ also separates $\widetilde{p}$ from $T_i$, $\widetilde{p}$ and $g(\widetilde{p})$ lie on the same side of $\widetilde{c}$. Moreover, $g(\widetilde{p})$ can not appear along the boundary of $\Delta_i$; if it did, by what we just showed it would lie on a stem side of $\Delta_i$, and then the corresponding stem will not be embedded, violating Definition \ref{def:stem system}.

If $g(\widetilde{p})$ is in the interior of $\Delta_i$,
then $g(\widetilde{\alpha}_i)$ and $g(\widetilde{\alpha}_{i+1})$ are arcs from $g(\widetilde{p})$ (inside $\Delta_i$) to the lifts $g(\widetilde{c}_i)$ and $g(\widetilde{c}_{i+1})$ in $\pi^{-1}\cC$ (outside $\Delta_i$-- indeed, by construction the interior of each secondary tulip is disjoint from each $\widetilde{c}_j$).
However, the endpoints of $g(\widetilde{c}_i)$ and $g(\widetilde{c}_{i+1})$ are not in the counterclockwise interval of $\partial_\infty\bH$ from $[g(\widetilde{c}_j)]^-=[\widetilde{c}_i]^-$ to $[g(\widetilde{c}_{j+1})]^+=[\widetilde{c}_{i+1}]^+$.
This implies that one of $\widetilde{\alpha}_i$ or $\widetilde{\alpha}_{i+1}$ intersects its $g$-translate, contradicting Definition \ref{def:stem system}.

If $g(\widetilde{p})$ is outside $\Delta_i$ but either to the left of $s_i$ or to the right of $s_{i+1}$ (as pictured in Figure~\ref{fig:rule out}), then
a stem side of $g(\Delta_j)$ would intersect a stem side of $\Delta_i$ while not forming any triangle with either base, violating \ref{it:stem triangles}. 

In the remaining case for the location of $g(\widetilde{p})$, we find that $\widetilde{p}$ would be interior to $g(\Delta_j)$, whence the argument is completed as executed when we dealt with the case $g(\widetilde{p})$ in the interior of $\Delta_i$ by reversing the roles of $i$ and $j$ and replacing $g$ with $g^{-1}$. 
\end{proof}

\begin{figure}
\vspace{.5cm}
\begin{lpic}{pics/secondaryNibsEmbed6(,8cm)}
\Large
\lbl[]{69,103;$F_i$}
\lbl[]{50,53;$\widetilde{p}$}
\lbl[]{45,30;$s_i$}
\lbl[]{95,30;$s_{i+1}$}
\lbl[]{9,95;$\widetilde{c}$}
\lbl[]{-10,130;$\widetilde{c}_{i+1}$}
\lbl[]{140,130;$\widetilde{c}_{i}$}
\end{lpic}
\caption{The `secondary tulip' $\Delta_i$, with its flower $F_i$ shaded, and base $\widetilde{c}$. If $g$ sends $T_j$ to $T_i$, taking base to base, then $g(\widetilde{p})$ can be ruled out from each of the regions below $\widetilde{c}$.}
\label{fig:rule out}
\end{figure}

\begin{proof}[Proof of Lemma~\ref{lem:flowers in same position}]
Fix a lift $\widetilde{c}$ of the base cuff, on the boundary of a lift $\widetilde{P}$ of $P$. As in \S\ref{subsec:type}, the tops of any flowers in $\widetilde{P}$, with base on $\widetilde{c}$, are in bijection with $\bZ[1/2]$.
Lemma~\ref{lem:flowers equal} implies that $\cF$ injects into $\pi_1S$-orbits of inner flowers with base at $\widetilde{c}$.
An inner flower at $\widetilde{c}$ is determined by the choice of the tops, a pair from $\bZ[1/2]$.
When the type is fixed, this pair is determined by the lesser of the two, and the action of $\Stab(\widetilde{c})$ on $\bZ[1/2]$ is translation by two.
Hence $|\cF|\le 4$ for each $d$.

Moreover, if $|\cF|>\lfloor 2/d \rfloor$, then there exist a pair of flowers $F_i,F_j\in \cF$, so that the $\pi_1S$-orbits contain representatives based along $\widetilde{c}$ whose tops alternate. 
But this implies that a pair of the main stem sides of a $\pi_1S$-translate of $\Delta_j$ would intersect a stem main side of $\Delta_i$ while not forming any triangle with either base, violating \ref{it:stem triangles}. \qedhere
\end{proof}

\begin{remark}
\label{rem:secondary nibs}
The construction of $\Omega_0$ above should remind the reader of exactly the construction pursued to define tulips and inner flowers.  
While this is the right idea, we acknowledge the technical annoyance of various subtle variations from our definitions of the floral invariants: the `secondary tulips' $\Delta_i$ above all have a vertex at $\widetilde{p}$. It is tempting to puncture the surface at $p$, making the $\Delta_i$ much closer to genuine nibs, in the sense of Przytycki. However, their `tops' still occur on lifts in $\pi^{-1}\cC$. 
Thus the `secondary tulips' $\Delta_i$ are intermediate between our tulips and Przytycki's nibs.
\end{remark}

\begin{remark} \label{rem: needformaximal} We demonstrate the necessity for using maximal tulips with the following example--see Figure \ref{fig:quadratic stems}. One imagines that each strand represents a stem at a point $p$ occurring somewhere to the far left. Each of the vertical separating curves should be interpreted as being in $\cC$, and one readily observes that many of these stems are prefixes of other stems in the picture. The stems are organized into linearly many (in $\chi$) parallel families, the $i^{th}$ of which contains $\approx i$ stems. If, on the other hand, one 
selects only one strand in each family 
(so that in the resulting system no stem is a prefix of another), we see that there are only linearly many in $\chi$.

\begin{figure}
\centering
\vspace{1cm}
\includegraphics[width=\linewidth]{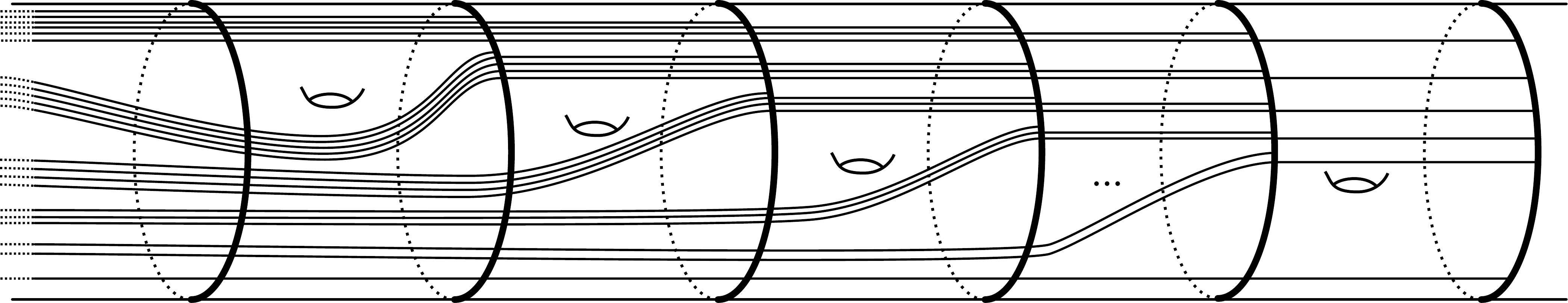}
\caption{If stems are allowed to be prefixes of each other in Definition~\ref{def:stem system}, there could be quadratically many in $|\chi|$, see Remark~\ref{rem: needformaximal}.}
\label{fig:quadratic stems}
\end{figure}

\end{remark}

We conclude this section with the lower bound in Proposition~\ref{prop:sigmag}: Suppose that $|\cC|=1$, and let $S_p$ be the surface obtained by removing $p$ and pinching $\cC$ to a pair of cusps; note that $|\chi(S_p)|=|\chi|+1$. Observe that  there are $2|\chi(S_p)|=
2(|\chi|+1)$ disjoint, pairwise non-homotopic arcs from $p$ to the cusps, which determines a stem system on $S$.  

\section{A technical lemma}
\label{sec:technical lemma}

The purpose of this section is to associate stem systems to the fiber over \emph{most} points (as measured by area) in the map from tulips to the hyperbolic surface $S$. 
Precisely, below we show that a point that lies in a family of inner flowers of tulips, all of the same refined type, gives rise to a pair of stem systems on $S$. 
We will deduce that the size of any fiber of the map from inner flowers to $S$ is at most $2\sigma(S)$, which is then controlled in \S\ref{sec:stems}.

\begin{lemma} \label{lem:tech2} 
Let $p\in S$, let $\widetilde{p}\in\bH$ be a lift of $p$, and let $\mathcal{A}_{\widetilde{p}}$ denote the set of stems  
at $\widetilde{p}$ of a collection of type-maximal tulips of common refined type, 
with $\widetilde{p}$ in their inner flower. Then the covering map induces $\mathcal{A}_{\widetilde{p}}\hookrightarrow\cA_1\sqcup\cA_2$, where each $\cA_i$ is a stem system based at $(p,\cC)$. 
\end{lemma}

\begin{figure}
\centering
\begin{lpic}{pics/stemIntersection(,10cm)}
\Large
\lbl[]{14,3;$\widetilde{c_-}$} 
\lbl[]{-3,94;$\widetilde{c_1}$}
\lbl[]{13,113;$\widetilde{d_1}$}
\lbl[]{75,3;$\widetilde{d_-}$}
\lbl[]{90,113;$\widetilde{c_2}$}
\lbl[]{103,95;$\widetilde{d_2}$}
\lbl[]{65,40;$t$}
\lbl[]{29,39;$s$}
\lbl[]{0,46;$\widetilde{c_0}$}
\lbl[]{13,23;$N$}
\lbl[]{79,23;$M$}
\end{lpic}
\caption{
A translate $g(t)$ that transversely intersects $s$, without intersecting $\widetilde{c_-}$.}
\label{fig:stemIntersection}
\end{figure}

\begin{proof}
Because $\widetilde{p}$ is in the inner flower of each of the tulips, by Proposition~\ref{prop:stems embed} each stem of $\cA_{\widetilde{p}}$ projects to a simple geodesic arc in $S$ from $p$ to $\cC$.
We first check that stems in $\cA_{\widetilde{p}}$ are pairwise non-homotopic \emph{rel $\cC$}; \ref{it:stem prefix} then follows immediately from the assumption that the stems of $\cA_{\widetilde{p}}$ are associated to type-maximal tulips. We will conclude by explaining why the image of $\cA_{\widetilde{p}}$ on $S$ can be partitioned into two pieces, each satisfying \ref{it:stem triangles}.

Suppose that $s,t\in\cA_{\widetilde{p}}$ are stems at $\widetilde{p}$ for tulips $N,M\to S$, whose inner flowers each contain $\widetilde{p}$, with tops $\widetilde{c_{1}},\widetilde{c_{2}}$ and $\widetilde{d_1},\widetilde{d_2}$, respectively, and (lifts of) main curves $\widetilde{\alpha_1}, \widetilde{\alpha_2}$ and $\widetilde{\beta_1}$, $\widetilde{\beta_2}$, respectively.
Suppose that $s$ and $t$ are homotopic \emph{rel $\cC$}.

Observe that the bases of $N$ and $M$ are the same curve $\widetilde{c_-}$.
This means that $\widetilde{\alpha_1}, \widetilde{\alpha_2}$ and $\widetilde{\beta_1}$, $\widetilde{\beta_2}$ are each consecutive in the linear order on $\Gamma_{\widetilde{c_-}}$.
However, because $\widetilde{p}$ is in both inner flowers, any geodesic in $\bH$ from $\widetilde{c_-}$ to $\partial_\infty\bH$ through $\widetilde{p}$ must separate both the tops of $N$ and those of $M$.
Therefore in the linear order $<$ on $\Gamma_{\widetilde{c_-}}$, both $\widetilde{\alpha_1}$ and $\widetilde{\beta_1}$ precede $\widetilde{\alpha_2}$ and $\widetilde{\beta_2}$; this is incompatible with each of the pairs $(\widetilde{\alpha_1},\widetilde{\alpha_2})$ and $(\widetilde{\beta_1},\widetilde{\beta_2})$ being consecutive.

The point $p$ lies in some pair of pants $P$, complementary to $\cP_\cC$. Because the stems arise from tulips of fixed refined type, each stem has a prefix which is a geodesic from $\widetilde{p}$ to a lift of a certain cuff $\gamma$ of $P$. Such geodesic prefixes come in two relative homotopy classes by Lemma~\ref{lem:homotopy classes}, which results in a partition of the image of $\cA_{\widetilde{p}}$ into $\cA_1$ and $\cA_2$, two collections of embedded geodesic arcs from $p$ to $\cC$, satisfying \ref{it:stem prefix}, whose elements all begin with a common prefix.

Let $\cA$ be one of these two subsets.
We now explain why pairs from $\cA$ satisfy \ref{it:stem triangles}.
Let $N,M\subset \bH$ be the tulips containing stems $s,t$ from the point $\widetilde{p}$, in the inner flower of both, with bases $\widetilde{c_-}$ and $\widetilde{d_-}$, and respective tops $\widetilde{c_1}$, $\widetilde{c_2}$ and $\widetilde{d_1}$, $\widetilde{d_2}$, labeled so that $\widetilde{c_1},\widetilde{c_-},\widetilde{c_2}$ and $\widetilde{c_1},\widetilde{d_-},\widetilde{c_2}$ are positively oriented in $\partial_\infty\bH$. 
By construction, $N$ and $M$ share the same flower base $\widetilde{c_0}$.
Without loss of generality, we may assume that $\widetilde{c_-},\widetilde{d_-},\widetilde{c_2}$ are positively oriented in $\partial_\infty\bH$.
Because $N$ and $M$ are of the same type, their tops alternate in the linear order on $\widetilde{c_0}$, i.e.~$\widetilde{c_1}\le \widetilde{d_1} < \widetilde{c_2}\le \widetilde{d_2}$.
Observe that this implies that the left main side of $M$ and the right main side of $N$ intersect.
See Figure~\ref{fig:stemIntersection}.

We proceed towards contradiction of \ref{it:stem triangles} by considering intersections between the $\pi_1S$-orbits of $s$ and $t$.  
If \ref{it:stem triangles} does not hold, there are points interior to $t$ that project to intersections of the images of $s$ and $t$, but are not vertices of a triangle on $S$ formed by the images of $s$ and $t$ and either base.
Consider such a point $O$, chosen to be closest to $\widetilde{p}$ along $t$. 
By assumption, there is $g\in\pi_1S$ so that $g(t)\pitchfork s$, and so that $g(t)\cap \widetilde{c_-}=\emptyset $.

Recall that, given a geodesic $L$ not through $\widetilde{p}$, we indicate the oriented interval in $\partial_\infty\bH$ on the other side of $L$ from $\widetilde{p}$ by $[L]$; we indicate the endpoints of $[L]$ by $[L]^\pm$. 
Observe that the image of the geodesic $\widetilde{d_-}$ under any covering transformation is disjoint from all other lifts in $\pi^{-1}\cP_\cC$. 
We will consider the location of $[g(\widetilde{d_-})]\subset \partial_\infty \bH$.
For simplicity, we will refer to the complementary component of $s\cup t\cup \widetilde{c_-}\cup\widetilde{d_-}$ cutting off a compact portion of $\widetilde{c_0}$ as the \emph{mountain}.

We make three claims:

\begin{enumerate}
\item $[g(\widetilde{d_-})]$ must lie in $[[\widetilde{c_2}]^-,[\widetilde{d_1}]^+]$, and $g(\widetilde{p})$ lies in the mountain.
\item $[g(\widetilde{c_-})]$ must lie in $[[\widetilde{c_2}]^-,[\widetilde{d_1}]^+]$.
\item $g(g(O))=O$.
\end{enumerate}

Of course, elements of $\pi_1S$ act as deck transformations and $\pi_1S$ is torsion-free. 
Thus claim (3) is a contradiction, and we may conclude that \ref{it:stem triangles} is satisfied. 

For claim (1): Because $g(t)\cap\widetilde{c_-}=\emptyset$ and because $g(\widetilde{p})$ is not contained in $[\widetilde{c_-}]$, the interval $[g(\widetilde{d_-})]$ is not contained in $[\widetilde{c_-}]$; because the main curves of $M$ intersect the projection of the base $\widetilde{d_-}$ once, $[g(\widetilde{d_-})]$ is not contained in $[\widetilde{d_-}]$ either.
If $[g(\widetilde{d_-})]$ lies in $[[\widetilde{c_-}]^+,[\widetilde{d_-}]^-]$, one finds a second intersection between the right side of $N$ and left side of $M$.
If $[g(\widetilde{d_-})]$ lies in $[[\widetilde{d_-}]^+,[\widetilde{c_2}]^-]$,
$g(t)\pitchfork t$.
If $[g(\widetilde{d_-})]$ lies in $[[\widetilde{d_1}]^+,[\widetilde{c_-}]^{-}]$, one first observes that $g(\widetilde{p})$ must lie in the mountain, since otherwise $g(t) \pitchfork t$. Moreover, $g(\widetilde{d_1})$ must also be in $[[\widetilde{d_1}]^+,[\widetilde{c_-}]^-]$ since the left main side of $g(M)$ can not cross the left main side of $M$, nor can it cross the right main side of $N$ since this would lead to two crossings between a main side of $M$ and a main side of $N$. But then $[g(\widetilde{c_{-}})]$ must also live in this same interval, since it occurs in $[[g(\widetilde{d_1})]^{+}, [g(\widetilde{d_{-}}]^{-}]$, 
which implies $g(s) \pitchfork s$. 
Therefore $[g(\widetilde{d_-})]$ is as claimed. 
Now $g(\widetilde p)$ is also as claimed:
$g(\widetilde{p})$ is to the right (as pictured in Figure~\ref{fig:stemIntersection}) of the geodesic containing $s$, to the left of $t$ (otherwise $g(t) \pitchfork t$), and on the same side of both $\widetilde{c_-}$ and $\widetilde{d_-}$ as $g(t)$.

For claim (2): 
If $[g(\widetilde{c_-})]$ lies in $[\widetilde{c_-}]$, the main curves of $g(N)$ intersect the projection of the base $\widetilde{c_-}$ more than once.
If $[g(\widetilde{c_-})]$ lies in $[[\widetilde{c_1}]^+,[\widetilde{c_-}]^-]$, $g(s)\pitchfork s$.
If $[g(\widetilde{c_-})]$ lies in $[[\widetilde{c_-}]^+,[\widetilde{c_2}]^-]$, then $(g(\widetilde{d_-}),g(\widetilde{d_1}))$ links with $(\widetilde{c_-},\widetilde{c_2})$, so the right side of $N$ intersects the left side of $M$ a second time.

Finally, for claim (3): By claims (1) and (2), the image $g(s\cup t)$ cuts off a quadrilateral with $s\cup t$. The choice of $O$ as closest to $\widetilde{p}$ along $t$ implies that the corners of this quadrilateral are $O$, $\widetilde{p}$, $g(O)$, and $g(\widetilde{p})$. The covering transformation $g$ must take the subarc $\widetilde{p}g(O)$ of $s$ to the fourth side $g(\widetilde{p})O$ of this quadrilateral, hence $g(g(O))=O$.
\qedhere

\end{proof}

\begin{remark}
\label{rem:better conclusion}
An earlier version of this manuscript claimed an improvement of Lemma~\ref{lem:tech2}, in that the partition of $\cA_{\widetilde{p}}$ was avoided, and no refined type assumption was imposed.
We are grateful to a careful referee who again spotted holes in that proof.
While we believe more involved bookkeeping would confirm the conclusion of that stronger lemma, we have chosen a streamlined exposition that avoids the necessary case-by-case analsis, at the expense, again, of worse constants.
\end{remark}

\section{Proofs of Theorems~\ref{thm:P1} and \ref{thm:P2}}
\label{sec:theorem proofs}

We collect the material from the previous sections to fill in the proofs of Theorems~\ref{thm:P1} and \ref{thm:P2}.

\begin{proof}[Proof of Theorem~\ref{thm:P1}]
Let $\cN_{(d,k)}$ be the collection of type-maximal tulips of refined type $(d,k)$. 

Consider the map
\[
\Phi:\bigsqcup_{N\in\cN_{(d,k)}} N^{\Delta} \looparrowright S~,
\]
where $N^{\Delta}$ is the inner flower of $N$.
For each $p\in S$, the fiber $\Phi^{-1}(p)$ forms a pair of
$(p,\cC)$ stem systems by Lemma~\ref{lem:tech2}, so is bounded in size by $2\sigma(S)$.
Considering the hyperbolic areas of the domain and codomain of $\Psi$, we find by Gau\ss-Bonnet that
\[
\sum_{N\in\cN_{(d,k)}} \Area\left(N^\Delta\right) \le \Area(S)\cdot 2\sigma(S) = 4\pi|\chi|\cdot \sigma(S)~.
\]
By Lemma~\ref{lem:inner flower area}, $\Area\left(N^\Delta\right)\to \pi$ as $\epsilon\to 0$, so it follows that $|\cN_{(d,k)}|\le 4|\chi|\sigma(S)$, whence $\cN=\sqcup \cN_{(d,k)}$ has size bounded by $84|\chi|\sigma(S)$ by Lemma~\ref{lem:refined types}.
\end{proof}

\begin{proof}[Proof of Theorem~\ref{thm:P2}]
Lemma~\ref{lem:accounting} demonstrates that each \emph{oriented} curve from $\Gamma\setminus\cC$ is the main side of a tulip that is maximal amongst its type, so each curve 
is the main side of a tulip that is maximal amongst its type \emph{at least two times}. Therefore $2|\Gamma\setminus\cC|$ is at most the number of sides of the set of tulips that are maximal amongst their type, and each tulip has two main sides. \qedhere
\end{proof}

\section{questions}
\label{sec:questions}

We conclude with some remaining questions. 
\bigskip

Of course, while our Main Theorem determines the correct growth rate of the maximum size of a $1$-system up to multiplicative constants, we still do not have the sort of exact answers that Przytycki obtains for arcs \cite{Przytycki}. So, first and foremost:

\begin{question} \label{q:exact}  
What is the exact maximum size of a $1$-system on $S$? Is it given by a polynomial expression in the genus and number of punctures?
\end{question}

Question \ref{q:exact} is already quite interesting for orientable closed surfaces, 
for which the only known answers have $g\le 2$:
On the torus the answer is $3$, and Malestein-Rivin-Theran determined the maximum size of a $1$-system in genus $2$ to be $12$ \cite{MRT}. 
On the other hand, in \cite{MRT} it is shown that the \textit{exact} answer for the size of a maximum $1$-system on the $n$-punctured torus is $3n$. Furthermore, letting $N(g,n)$ denote the maximum size of a $1$-system on a surface of genus $g \ge 2$ and with $n \ge 0$ punctures, they show that 
\[ (2g+1)(n+1) \leq N(g,n) \leq N(g,0) + (2g+1)\, n. \]
In particular, if one imagines each lattice point in the first quadrant of the $(g,n)$-plane as corresponding to one such surface, the maximum size of a $1$-system grows only linearly along vertical lines. On the other hand, examples demonstrate that it must grow quadratically along horizontal lines, and in fact along any line along which $g \rightarrow \infty$. In any case, because the maximum size evidently behaves differently with the addition of punctures as opposed to genus, we frame Question~\ref{q:exact} in terms of $g$ and $n$ instead of merely the single variable $\chi$.

Recently, Zhai \cite{Zhai} produced a $1$-system on the closed genus $3$ surface with $33$ curves and proved that it is maximal with respect to inclusion. 
Nonetheless, $N(3,0)$ remains unknown. 

\begin{question} \label{q:genus 3} What is the maximum size of a $1$-system on a closed orientable surface of genus $3$?
\end{question}

We note that our Main Theorem produces the upper bound $33606$ when $g=3$, evidently far too large. For example, the exponential-in-genus bound $(g-1)(4^g-1)$ produced by Malestein-Rivin-Theran \cite{MRT} is bigger than our quadratic-in-genus bound starting only when $g=8$; when $g=3$, their bound yields $126$.

There are several obvious routes forward for sharpening our bounds: 

\begin{question}[Improve Lemma~\ref{lem:accounting}] \label{q:maximal} Is every curve in a $1$-system a main curve of a maximal tulip?
\end{question}

\begin{question}[Improve Lemma~\ref{lem:tech2}]
\label{q:refined}
Do the stems associated to points in inner flowers of a collection of maximal tulips form a stem system (in the sense of Definition~\ref{def:stem system})? 
\end{question}

\begin{question}[Improve Proposition~\ref{prop:stembounds}]
\label{q:sigma} What is the value of $\sigma(S)$ exactly?
\end{question}

If the answer to Question~\ref{q:maximal} is yes, one can delete a factor of $7$ from the righthand side of \eqref{eq:main inequality};
if the answer to Question~\ref{q:refined} is yes, one can delete a factor of $6$;
if the answer to Question~\ref{q:sigma} is $2|\chi|+2$, one can delete a factor of $12$.
In the hopeful scenario of all three answers, our bound in genus $3$ would improve to $6+4\cdot4\cdot5=86$, an improvement on the Malestein-Rivin-Theran bound.

Even more ambitiously, it is entirely conceivable (likely even?) that each curve in a 1-system is the main curve of \emph{many} maximal tulips, since a typical element of $\Gamma$ will intersect many curves in $\cC$, and many other curves in $\Gamma$. More intersections yields more complicated looking $1$-systems, which should increase the number of maximal tulips.
With that in mind, perhaps each curve is one of the main curves of \emph{at least four} maximal tulips-- this is after all the situation for arcs. Following this stronger bound through the proof of \eqref{eq:main inequality}, together with the optimistic answers to the above questions, would yield the dramatically stronger bound $6+\frac12\cdot4\cdot4\cdot 5=46$ in genus $3$, evidently not far from Zhai's lower bound. Of course, such improvements are interesting in higher genus as well, as they would decrease the coefficient of $|\chi|^2$.

As mentioned in the introduction, to the best of our knowledge, all previous attempts at proving upper bounds for the size of a $1$-system relied on treating Przytycki's bounds for arcs as a black box. The idea in all of these attempts is to start with a system of curves and somehow derive a system of arcs from it, perhaps by puncturing the surface somewhere and ``dragging'' each curve towards it. The key difficulty here is that some of the curves can obstruct this puncture from the other curves, in such a fashion that the only way one can drag all of the curves into it is to force additional intersections between curves. One might try to mitigate this by adding more than one puncture, but then of course Przytycki's bound for arcs depends on the Euler characteristic of the underlying surface and will weaken as one adds punctures.

In any case, we will informally refer to this type of argument as \textit{arcification}, and we ask whether there exists a way to obtain quadratic-in-$\chi$ upper bounds with an arcification-style proof:

\begin{question} \label{q:arcification} Is there a way to achieve sharp upper bounds on the size of a $1$-system via arcification? If not, is there a combinatorial characterization of those $1$-systems that can be arcified? 
\end{question}

As mentioned several times above, the proof of our Main Theorem is complicated by the possibility $\cC\subsetneq \cP_\cC$; in other words, $\Gamma$ may not contain a pants decomposition. We still have no idea what $1$-systems of maximum cardinality really look like, in particular: 

\begin{question}
\label{q:pants} 
Does every $1$-system of maximum cardinality contain
a pants decomposition? 
\end{question}

Amazingly, we even do not know if every maximal $1$-system fills the surface.
Here we mean `maximal' with respect to inclusion among $1$-systems-- this is called \emph{saturated} in \cite{Zhai}-- and `filling' in the sense that cutting along them produces a disjoint union of topological disks and once-punctured disks. (It is known that the only other possibility is that there might be some annular components-- see \cite{MRT}.) 

\begin{question} \label{q:fill} Does every maximal $1$-system fill the surface on which it lives? 
\end{question}

We imagine that the answer to Question \ref{q:fill} must be ``yes''. Nonetheless, here is how one might imagine it could be ``no'': 
It seems conceivable that a $1$-system $\Gamma$ has the property that if one cuts along a curve $\gamma\in \Gamma$, then any arc with one endpoint on each of the two new boundary components has the property that it intersects at least one curve in $\Gamma$ at least two times. 

Finally, we observe that all known constructions of large $1$-systems are of an inductive flavor: one begins with a $1$-system on a surface of small Euler characteristic. One then adds $1$-handles while also carefully adding in additional curves in such a way that the addition of the $n^{th}$ $1$-handle allows for at least (roughly) $n$ new curves, guaranteeing that the size of the $1$-system is always at least quadratic in the Euler characteristic of the underlying surface. But perhaps there are very large $1$-systems that can not be built up in this fashion, and so we ask: 

\begin{question} \label{q:inductive} Is there a maximum cardinality $1$-system on a closed orientable surface that is not ``inductively constructed'' as described above? Perhaps this can be interpreted as asking whether there is such a $1$-system $\Gamma$ so that cutting along any non-separating curve $\beta$ (which may or may not be in $\Gamma$), capping off the resulting boundary components, and deleting any curve that intersected $\beta$, yields a $1$-system on a genus $g-1$ surface that is far from maximum size. 
\end{question}

\bigskip
\bibliography{references}

@article {Przytycki,
    AUTHOR = {Przytycki, Piotr},
     TITLE = {Arcs intersecting at most once},
   JOURNAL = {Geometric and Functional Analysis},
  FJOURNAL = {Geometric and Functional Analysis},
    VOLUME = {25},
      YEAR = {2015},
    NUMBER = {2},
     PAGES = {658--670},
      ISSN = {1016-443X,1420-8970},
   MRCLASS = {57M50 (05B99)},
  MRNUMBER = {3334237},
MRREVIEWER = {Aleksander\ Malni\v c},
       DOI = {10.1007/s00039-015-0320-0},
       URL = {https://doi.org/10.1007/s00039-015-0320-0},
}

@book{Buser,
  title={Geometry and spectra of compact Riemann surfaces},
  author={Buser, Peter},
  year={2010},
  publisher={Springer Science \& Business Media}
}

@article {JMM,
    AUTHOR = {Juvan, M. and Malni\v{c}, A. and Mohar, B.},
     TITLE = {Systems of curves on surfaces},
   JOURNAL = {Journal of Combinatorial Theory. Series B},
  FJOURNAL = {Journal of Combinatorial Theory. Series B},
    VOLUME = {68},
      YEAR = {1996},
    NUMBER = {1},
     PAGES = {7--22},
      ISSN = {0095-8956,1096-0902},
   MRCLASS = {57M15 (05C10)},
  MRNUMBER = {1405702},
       DOI = {10.1006/jctb.1996.0053},
       URL = {https://doi.org/10.1006/jctb.1996.0053},
}

@article {MRT,
    AUTHOR = {Malestein, Justin and Rivin, Igor and Theran, Louis},
     TITLE = {Topological designs},
   JOURNAL = {Geometriae Dedicata},
  FJOURNAL = {Geometriae Dedicata},
    VOLUME = {168},
      YEAR = {2014},
     PAGES = {221--233},
      ISSN = {0046-5755,1572-9168},
   MRCLASS = {57M50 (05B99 57M15)},
  MRNUMBER = {3158040},
MRREVIEWER = {Dai\ Tamaki},
       DOI = {10.1007/s10711-012-9827-9},
       URL = {https://doi.org/10.1007/s10711-012-9827-9},
}

@article {ABG,
    AUTHOR = {Aougab, Tarik and Biringer, Ian and Gaster, Jonah},
     TITLE = {Packing curves on surfaces with few intersections},
   JOURNAL = {International Mathematics Research Notices (IMRN)},
  FJOURNAL = {International Mathematics Research Notices. IMRN},
      YEAR = {2019},
     PAGES = {5205--5217},
      ISSN = {1073-7928,1687-0247},
   MRCLASS = {57K32},
  MRNUMBER = {4001026},
MRREVIEWER = {Hongbin\ Sun},
       DOI = {10.1093/imrn/rnx270},
       URL = {https://doi.org/10.1093/imrn/rnx270},
}

@article {Greene1,
    AUTHOR = {Greene, Joshua Evan},
     TITLE = {On loops intersecting at most once},
   JOURNAL = {eometric and Functional Analysis},
  FJOURNAL = {Geometric and Functional Analysis},
    VOLUME = {29},
      YEAR = {2019},
    NUMBER = {6},
     PAGES = {1828--1843},
      ISSN = {1016-443X,1420-8970},
   MRCLASS = {57M15 (05C62 05D40)},
  MRNUMBER = {4034921},
MRREVIEWER = {Lorenzo\ Traldi},
       DOI = {10.1007/s00039-019-00517-0},
       URL = {https://doi.org/10.1007/s00039-019-00517-0},
}

@book {Farb-Margalit,
    AUTHOR = {Farb, Benson and Margalit, Dan},
     TITLE = {A primer on mapping class groups},
    SERIES = {Princeton Mathematical Series},
    VOLUME = {49},
 PUBLISHER = {Princeton University Press, Princeton, NJ},
      YEAR = {2012},
     PAGES = {xiv+472},
      ISBN = {978-0-691-14794-9},
   MRCLASS = {57M50 (20F36 20F65 57M07 57N05)},
  MRNUMBER = {2850125},
MRREVIEWER = {Stephen\ P.\ Humphries},
}

@article {Zhai, 
    AUTHOR = {Zhai, Zhaoshen}, 
    TITLE = {A saturated 1-system of curves on the surface of genus 3}, 
    JOURNAL = {Geometriae Dedicata}, 
    VOLUME = {219}, 
    NUMBER = {49}, 
    YEAR = {2025}, 
}

@book{k3,
  title={K3: a new problem list in low-dimensional topology},
  author={Baykur, Refik {\.I}nan{\c{c}} and Kirby, Robion C and Ruberman, Daniel},
  volume={295},
  year={2026},
  publisher={American Mathematical Society}
}

\bibliographystyle{alpha}

\end{document}